\documentclass[11pt]{article}
\usepackage{bm}
\usepackage{fullpage,graphicx,psfrag,amsmath,amsfonts,verbatim}
\usepackage{hyperref}
\usepackage{enumerate}
\usepackage{subfig}
\usepackage{multirow}
\usepackage{algorithm}
\usepackage{algpseudocode}

\title{{Relax-Tighten-Round Algorithm for Optimal Placement and Control of Valves and Chlorine Boosters in Water Networks}}

\author{Filippo Pecci$^1$, Ivan Stoianov$^1$, Avi Ostfeld$^2$}
\date{%
    $^1$InfraSense Labs, Department of Civil and Environmental Engineering, Imperial College London, London, SW7 2BB, UK\\%
    $^2$Faculty of Civil and Environmental Engineering Technion - Israel Institute of Technology Haifa 32000, Israel\\[2ex]%
    March 2021
}

\begin{document}
\maketitle
\begin{abstract}
In this paper, a new mixed integer nonlinear programming formulation is proposed for optimally {placing} and operating pressure reducing valves and chlorine booster stations in water distribution networks. The objective is the minimization of average zone pressure, while penalizing deviations from a target chlorine concentration. 
{We propose a relax-tighten-round algorithm based on tightened polyhedral relaxations and a rounding scheme to compute feasible solutions, with bounds on their optimality gaps. This is because off-the-shelf global optimization solvers failed to compute feasible solutions for the considered non-convex mixed integer nonlinear program.} The implemented algorithm is evaluated using three benchmarking water networks, and they are shown to {outperform} off-the-shelf solvers, for these case studies. The proposed heuristic has enabled the computation of good quality feasible solutions in most instances, {with bounds on the optimality gaps that are comparable to the order of uncertainty observed in operational water network models.}
\end{abstract}

\section{Introduction}
{The main} operational objectives for water utilities include {the} reduction of water leaks, and management of drinking-water quality. Leakage reduction is achieved by controlling average zone pressure (AZP) within water distribution networks (WDNs), while satisfying minimum service requirements \cite{Wright2015}. Pressure control schemes are implemented through pressure reducing valves (PRVs), which reduce pressure at their downstream node. The problem of optimal placement and operation of PRVs in WDNs has been formulated in previous literature, and solved using both evolutionary algorithms \cite{Araujo2006,Nicolini2009} and mathematical optimization methods~\cite{Eck2012,PecciGlobalOpt}. 

Monitoring and control of disinfectant residuals {in} drinking water distribution networks {is critical to maintain the water quality and eliminate the risks of contamination with pathogens such as bacteria and viruses in distribution ~\cite{Aisopou2014AnalyticalConditions,Sakomoto2020ManagingUganda}. This is particularly critical during the current COVID-19 pandemic as leaking sewage from sewer networks could allow potentially harmful contaminants into drinking water networks \cite{Quilliam2020COVID-19:Faeces}.
In order to deactivate any pathogens that might exist in distribution networks}, disinfectant is typically added at water sources (e.g. water treatment plants), with chlorine being a commonly used water disinfectant. Because chlorine is reactive, it is depleted over time as it travels {across the pipe networks}, causing a reduction in the ability to prevent microbial contamination. Water utilities aim to maintain a target chlorine concentration, which is sufficient to safeguard public health, while avoiding excessive chlorination, resulting in taste and odor problems, as well as {the growth} of {disinfection} {by-products}. 
{In addition, the objective is to maintain optimal and constant chlorine concentrations, as variations in chlorine concentration are perceived as water quality problems by customers.}
Chlorine booster stations {are used} to deal with this challenge \cite{Boccelli1998,Propato2004a}. 
Using booster chlorination, disinfectant is re-applied at selected locations within the network, leading to a more uniform spatial and temporal distribution of chlorine residuals. Previous literature has modeled the operation of booster stations assuming known flow velocities across network pipes - see as examples \cite{Boccelli1998,Propato2004a}. However, this can lead to sub-optimal design and operation of WDNs. In fact, in order to mitigate disinfectant decay reactions, network operators should aim to reduce travel time from water sources to demand nodes. This may {result} in sub-optimal pressure management schemes, where minimum pressure constraints are not satisfied, as observed in \cite{Kang2010}. Therefore, 
we consider the joint optimization of hydraulic pressure and flows, together with chlorine {residual} concentrations in WDNs. 

We investigate the problem of minimizing average zone pressure, while penalizing {deviations} from chlorine target concentrations, and satisfying regulatory constraints on pressure and chlorine concentration levels. \cite{Ostfeld2005} and \cite{Kang2010} implemented genetic algorithms to solve problems of optimal operation of WDNs, where optimization unknowns include network flows and chlorine concentrations, while locations of PRVs and chlorine booster stations are fixed. However, pressure reducing valves and booster stations should be optimally placed for a more effective pressure control and {management of chlorine residual concentrations}. 

In this manuscript, we propose a new mathematical framework for the optimal placement and operation of pressure reducing valves and chlorine booster stations.
The considered objective is the minimization of average zone pressure, while penalizing {deviations} from target chlorine concentrations at demand nodes. 
The transport of chlorine through each pipe is modeled by a one dimensional first-order advection PDE \cite{Rossman1996}, where flow velocity corresponds to the one-dimensional velocity field, and a linear function is used to represent chlorine decay \cite{Hallam2002TheSystems}. We implement an implicit upwind scheme to discretize the considered PDE. Optimization constraints include quadratic equations modelling head loss due to pipe friction \cite{Eck2015,PecciQuadApp}, and bilinear terms due to the presence of unknown flow velocities within the discretized advection PDE. In addition, binary variables are used to model the direction of flow across pipes, and the placement of valves and booster stations. The resulting optimization problem is a non-convex mixed integer nonlinear program. 

In comparison to previous literature \cite{Ostfeld2005,Kang2010}, which relied on genetic algorithms, we investigate the application of mathematical optimization methods to compute feasible solutions for the considered problem, with guaranteed bounds on their optimality gaps.
{We propose a relax-tighen-round (RTR) algorithm based on polyhedral relaxations of the non-convex terms, an optimization-based-bound-tightening scheme, and a rounding heuristic. The developed RTR algorithm computes a feasible solution for the considered non-convex MINLP, with bounds on its optimality gap.} {In comparison, we show that off-the-shelf global optimization solvers failed to generate feasible solutions for the considered problem. }The performance of the { RTR algorithm} is investigated using multiple problem instances for different WDN case studies.

\section{Problem formulation}
We formulate the problem of optimal placement and operation of pressure reducing valves and chlorine booster stations, with the objective of minimizing average zone pressure, while penalizing deviations from target chlorine concentrations at demand nodes. 
A WDN with $n_n$ demand nodes, $n_0$ source nodes (e.g. water sources, water treatment plants), and $n_p$ links is modelled as a directed graph with $n_n+n_0$ vertices and $n_p$ edges. Define $\mathcal{P}:=\{1,\ldots,n_p\}$ and $\mathcal{N}:=\{1,\ldots,n_n\}$, $\mathcal{N}^0:=\{1,\ldots,n_0\}$.
Given a node $i \in \mathcal{N}$, let $I^{\text{in}}_i$ and $I^{\text{out}}_i$ be the index sets corresponding to links with assigned direction entering and leaving the node, respectively.
We consider network operation within a discretized time interval $\mathcal{T} = {1,\ldots,n_t}$.
{The objective of this study is to minimize average zone pressure in water distribution networks, while penalizing deviation from target chlorine concentrations. Average Zone Pressure (AZP) is defined as the following weighted sum of nodal pressures \cite{Wright2015}:
\begin{equation}
\label{eq:AZP}
    \sum_{k \in \mathcal{T}}\sum_{i \in \mathcal{N}}\omega_i(h_{i,k}-h^{\text{elev}}_i)
\end{equation}
where $h_{i,k}$ is the unknown hydraulic head at node $i \in \mathcal{N}$ and time $k \in \mathcal{T}$, while $\bm{h}^{\text{elev}} \in \mathbb{R}^{n_n}$ is the vector of known nodal elevations. Weights are defined as follows:
\begin{equation}
    \omega_i:= \frac{\sum_{l \in I^{\text{in}}_i \cup I^{\text{out}}_i}L_l }{n_t\sum_{j \in \mathcal{N}}\sum_{l \in I^{\text{in}}_j \cup I^{\text{out}}_j}L_l }, \quad i\in \mathcal{N}
\end{equation}
Let $\bm{c}^* \in \mathbb{R}^{n_n}$ be a vector of target chlorine concentration at network nodes. Moreover, set
\begin{equation}
    \hat{d}_{i,k}=\frac{d_{i,k}}{\sum_{k \in \mathcal{T}}\sum_{j \in \mathcal{N}}d_{j,k}}, \quad i \in \mathcal{N},\; k \in \mathcal{T}
\end{equation}
where $d_{i,k}$ is the known demand at node $i \in \mathcal{N}$ and time $k \in \mathcal{T}$. Denote by $c_{i,k}$ the unknown chlorine concentration at node $i \in \mathcal{N} \cup \mathcal{N}^0$ and time $k \in \mathcal{T}$. We define the Average Target Deviation (ATD) as 
\begin{equation}
\label{eq:ATD}
    \sum_{k \in \mathcal{T}}\sum_{i \in \mathcal{N}}\hat{d}_{i,k}|c_{i,k}-c^*_i|
\end{equation}
The formula for ATD can be reformulated as a linear function by introducing auxiliary variables $\mu_{i,k}$, which satisfy the following linear constraints:
\begin{subequations}
\label{eq:ATD_auxvar}
\begin{align}
    &c_{i,k}-c_i^* \leq \mu_{i,k},\quad i \in \mathcal{N},\; k\in\mathcal{T} \\
    &-c_{i,k}+c_i^* \leq \mu_{i,k},\quad i \in \mathcal{N},\; k\in\mathcal{T}.
    \end{align}
\end{subequations}
The objective function to be minimized is written as:
\begin{equation}
\label{eq:objfun}
        \sum_{k \in \mathcal{T}}\sum_{i \in \mathcal{N}}\omega_ih_{i,k} + \sum_{k \in \mathcal{T}}\sum_{i \in \mathcal{N}}\hat{d}_{i,k}\mu_{i,k}
\end{equation}
Since the considered problem aims to optimize both hydraulic pressure and water quality, its formulation is based on hydraulics and water quality modelling. }
\subsection{Hydraulic variables and constraints}
First, we introduce optimization variables and constraints related to network hydraulic properties. Source nodes are assumed to have known hydraulic heads $h^0_{i,k}$, $i \in \mathcal{N}^0$, $k \in \mathcal{T}$. We denote by $q_{l,k}$ the unknown flow in link $l \in \mathcal{P}$ at time $k \in \mathcal{T}$. The unknown frictional head loss across link $l$ at time $k$ is denoted by $\theta_{l,k}$.
Pressure control valves reduce pressure at their downstream node, introducing additional head losses, which are represented by variable $\eta_{l,k}$, $l\in\mathcal{P}$, $k \in\mathcal{T}$. Vector of binary variables $\bm{v} \in \{0,1\}^{2n_p}$ models the placement of control valves. We have:
\begin{equation}
    v_l = \begin{cases} 1 &\text{a valve is placed on link $l$ in the positive flow direction} \\
    0 & \text{otherwise},
    \end{cases}
\end{equation}
and
\begin{equation}
    v_{n_p+l} = \begin{cases} 1 &\text{a valve is placed on link $l$ in the negative flow direction} \\
    0 & \text{otherwise}.
    \end{cases}
\end{equation}
These binary variables are subject to the following physical and economical constraints:
\begin{subequations}
\label{eq:binary_cons}
\begin{align}
&v_l+v_{n_p+l} \leq 1, \quad l \in \mathcal{P} \\
&\sum_{l\in \mathcal{P}}(v_l+v_{n_p+l}) = n_v
\end{align}
\end{subequations}
The following constraints formulate energy and mass conservation laws, and the placement of pressure reducing valves on network links:
\begin{subequations}
\label{eq:hyd_lincon}
\begin{align}
&h_{i_1,k} - h_{i_2,k} = \theta_{l,k} + \eta_{l,k}, \quad   i_1 \xrightarrow{l} i_2, l \in \mathcal{P}, i_1\in \mathcal{N},i_2 \in \mathcal{N}, k\in \mathcal{T}\\
&h^0_{i_1,k} - h_{i_2,k} = \theta_{l,k} + \eta_{l,k}, \quad   i_1 \xrightarrow{l} i_2, l \in \mathcal{P}, i_1\in \mathcal{N}^0,i_2 \in \mathcal{N}, k\in \mathcal{T}\\
&\sum_{l \in I^{\text{in}}_i}q_{l,k} - \sum_{l \in \in I^{\text{out}}_i}q_{l,k}= d_{i,k}, \quad i \in \mathcal{N}, k \in \mathcal{T}\label{eq:hyd_lin_con_mass_bal}\\
&\eta_{l,k} -\eta^{\max}_{l,k}v_l \leq 0, \quad l \in \mathcal{P}, \; k \in \mathcal{T}\\
&-\eta_{l,k}+\eta^{\min}_{l,k}v_{n_p+l} \leq 0, \quad l \in \mathcal{P}, \;k\in \mathcal{T}\\
&-q_{l,k} - q^{\min}_{l,k}v_l \leq -q_{l,k}^{\min}, \quad l \in \mathcal{P}, \;k \in \mathcal{T}\\
&q_{l,k} + q^{\max}_{l,k}v_{n_p+l} \leq q_{l,k}^{\max}, \quad l \in \mathcal{P}, \;k \in \mathcal{T}
\end{align}
\end{subequations}
In order to model the transport of chlorine constituent, it is required to explicitly {consider the flow direction} across network links as a decision variable. Therefore, we introduce auxiliary variables $q^+_{l,k}$, $q^-_{l,k}$, $\theta^+_{l,k}$, $\theta^-_{l,k}$, $s_{l,k}$, and binary variable $z_{l,k} \in \{0,1\}$ such that
\begin{subequations}
\label{eq:hyd_aux_con}
\begin{align}
& q_{l,k} = q^+_{l,k}-q^-_{l,k}, \quad l\in \mathcal{P},\; k\in \mathcal{T}\\
& s_{l,k} = q^+_{l,k}+q^-_{l,k}, \quad  l \in \mathcal{P}, \; k\in \mathcal{T} \\
& \theta_{l,k} = \theta^+_{l,k} - \theta^-_{l,k},  \quad l\in \mathcal{P},\; k\in \mathcal{T}\\
& 0 \leq q^+_{l,k}\leq (q^{+})^{\max}_{l,k}z_{l,k}, \quad l\in \mathcal{P},\; k \in \mathcal{T} \\
& 0 \leq q^-_{l,k}\leq (q^-)^{\max}_{l,k}(1-z_{l,k}), \quad l\in \mathcal{P},\; k \in \mathcal{T} \\
& 0 \leq \theta^+_{l,k}\leq (\theta^{+})^{\max}_{l,k}z_{l,k}, \quad l\in \mathcal{P},\; k \in \mathcal{T} \\
&0 \leq \theta^-_{l,k}\leq (\theta^-)^{\max}_{l,k}(1-z_{l,k}), \quad l\in \mathcal{P},\; k \in \mathcal{T} 
\end{align}
\end{subequations}
Frictional head losses are often represented by either {the} Hazen-Williams (H-W) or {the} Darcy-Weisbach (D-W) {equations}~\cite{DAmbrosio2015}. Since both formulae involve non-smooth terms, quadratic approximations have been proposed and used in previous literature \cite{Eck2015,PecciQuadApp}. Let $\bm{a} \in \mathbb{R}^{n_p}$ and $\bm{b} \in \mathbb{R}^{n_p}$ be vector of coefficients of these approximations. We enforce the following constraints on variables $\theta^+_{l,k}$ and $\theta^-_{l,k}$:
\begin{subequations}
\label{eq:hyd_nonconvex}
\begin{align}
    \theta^+_{l,k} = a_l(q^+_{l,k})^2+ b_lq^+_{l,k}, \quad l\in \mathcal{P},\; k\in \mathcal{T}\label{eq:hyd_nonconvexa}\\
    \theta^-_{l,k} = a_l(q^-_{l,k})^2+ b_lq^-_{l,k}, \quad l\in \mathcal{P},\; k\in \mathcal{T}\label{eq:hyd_nonconvexb}.
\end{align}
\end{subequations}
{Constraints $z_{l,k} \in \{0,1\}$, \eqref{eq:hyd_aux_con}, and \eqref{eq:hyd_nonconvex} are equivalent to the non-linear equations:
\begin{subequations}
\label{eq:hyd_nlp}
\begin{align}
    &s_{l,k} = |q_{l,k}|, \quad l\in \mathcal{P},\; k\in \mathcal{T}\\
    &q^+_{l,k} = \max(q_{l,k},0), \quad l\in \mathcal{P},\; k\in \mathcal{T}\\
    &q^-_{l,k} = -\min(q_{l,k},0), \quad l\in \mathcal{P},\; k\in \mathcal{T}\\
    &\theta^+_{l,k} = \max(\theta_{l,k},0), \quad l\in \mathcal{P},\; k\in \mathcal{T}\\
    &\theta^-_{l,k} = -\min(\theta_{l,k},0), \quad l\in \mathcal{P},\; k\in \mathcal{T},\\
    &z_{l,k} = \frac{1+\text{sign}(q_{l,k})}{2}, \quad l \in \mathcal{P},\; k \in \mathcal{T},
\end{align}
\end{subequations}
and
\begin{equation}
\label{eq:hyd_nlpa}
        \theta_{l,k} = a_l|q_{l,k}|q_{l,k}+ b_lq_{l,k}, \quad l\in \mathcal{P},\; k\in \mathcal{T},
\end{equation}
where $\text{sign}(q_{l,k})=1$ if $q_{l,k}>0$, and $\text{sign}(q_{l,k})=-1$ otherwise.}
Lower and upper bounds on hydraulic variables are given by
\begin{subequations}
\label{eq:hyd_bounds1}
\begin{align}
    &q^{\min}_{l,k} \leq q_{l,k} \leq q^{\max}_{l,k}, \quad l\in \mathcal{P},\; k \in \mathcal{T} \\
    &h^{\min}_{i,k} \leq h_{i,k} \leq h^{\max}_{i,k}, \quad i\in \mathcal{N},\; k \in \mathcal{T} \\
    &\eta^{\min}_{l,k} \leq \eta_{l,k} \leq \eta^{\max}_{l,k}, \quad l\in \mathcal{P},\; k \in \mathcal{T}\\
   &\theta^{\min}_{l,k} \leq \theta_{l,k} \leq \theta^{\max}_{l,k}, \quad l\in \mathcal{P},\; k \in \mathcal{T}
\end{align}
\end{subequations}
and
\begin{subequations}
\label{eq:hyd_bounds2}
\begin{align}
    &(q^+)^{\min}_{l,k} \leq q^+_{l,k} \leq (q^+)^{\max}_{l,k}, \quad l\in \mathcal{P},\; k \in \mathcal{T} \\
    &(q^-)^{\min}_{l,k} \leq q^-_{l,k} \leq (q^-)^{\max}_{l,k}, \quad l\in \mathcal{P},\; k \in \mathcal{T} \\
    &(\theta^+)^{\min}_{l,k} \leq \theta^+_{l,k} \leq (\theta^+)^{\max}_{l,k}, \quad l\in \mathcal{P},\; k \in \mathcal{T} \\
    &(\theta^-)^{\min}_{l,k} \leq \theta^-_{l,k} \leq (\theta^-)^{\max}_{l,k}, \quad l\in \mathcal{P},\; k \in \mathcal{T}\\
    &s^{\min}_{l,k} \leq s_{l,k} \leq s^{\max}_{l,k}, \quad l\in \mathcal{P},\; k \in \mathcal{T}.
\end{align}
\end{subequations}
\subsection{Water quality variables and constraints}
Next, we describe variables and constraints associated with water quality. Let $c^{\max}_{i,k}$ be the maximum allowed chlorine concentration at network node $i \in \mathcal{N}\cup \mathcal{N}^0$ and time $k \in \mathcal{T}$. The evolution of chlorine concentration along a given link $l \in \mathcal{P}$ is governed by a PDE modelling advective transport of constituent with first order decay \cite{Rossman1996}. We implement an Eulerian implicit upwind discretization scheme, where backward differences are used to approximate both temporal and spatial derivatives \cite{Islam1999}. For each link $l \in\mathcal{P}$, we introduce a space discretization $j\Delta x_l$, $j \in \{0,\ldots,J_l\}$, with $\Delta x_l =\frac{L_l}{J_l}$, where $L_l$ is the length of link $l$. We denote by $r_{j,l,k}$ the chlorine concentration at $j\Delta x_l$ and time $k$, for all $j \in \{0,\ldots,J_l\}$ and $k \in \{0,\ldots,n_t\}$. We also have auxiliary variables $w_{j,l,k}$ such that:
\begin{equation}
\label{eq:wq_nonconvex}
    w_{j,l,k} = s_{l,k}r_{j,l,k}, \quad j=0,\ldots,J_l, \; l \in \mathcal{P}, \;  k \in \mathcal{T},
\end{equation}
For all $j \in \{1,\ldots,J_l\}$, $l \in \mathcal{P}$, and $k \in \mathcal{T}$, the discretized PDE yields:
\begin{equation}
\label{eq:upwind_disc}
    \begin{split}
   (1+\alpha_l \Delta t)r_{j,l,k} - r_{j,l,k-1} + \gamma_l(w_{j,l,k}-w_{j-1,l,k}) = 0,
    \end{split}
\end{equation}
where $\gamma_l = (4\Delta t)/(10^3\pi D^2_l \Delta x_l)$, with $L_l$ and $D_l$ length and diameter of link $l$, respectively, and $\alpha_l>0$ is the first order decay coefficient associated with pipe $l$ \cite{Hallam2002TheSystems}. Initial concentrations in pipes are defined as 
\begin{equation}
    r_{j,l,0} = c^0_{i_2}, \quad \forall j \in \{0,\ldots,J_l\}, \; \forall l \in \mathcal{P}, \quad i_1 \xrightarrow{l}{i_2}
\end{equation}
with given initial concentration $c^0_i$  at node $i \in \mathcal{N}\cup\mathcal{N}^0$. Furthermore, $r_{0,l,k}$ is assumed to be equal to the concentration of the upstream node, depending on the flow direction:
\begin{subequations}
    \label{eq:boundary_disc}
\begin{align}
&r_{0,l,k} - c_{i_1,k} + c^{\max}_{i_2}z_{l,k} \leq c^{\max}_{i_2}\\
&-r_{0,l,k} + c_{i_1,k} + c^{\max}_{i_1}z_{l,k} \leq   c^{\max}_{i_1}\\
&r_{0,l,k} - c_{i_2,k} - c^{\max}_{i_1}z_{l,k} \leq 0 \\
&-r_{0,l,k} + c_{i_2,k} -c^{\max}_{i_2}z_{l,k} \leq 0 ,
\end{align}
\end{subequations}
Our problem formulation considers as free decision variables concentrations at source nodes $c_{i,k}$, $i \in \mathcal{N}^0$. Moreover, let $\bm{v}^b \in \{0,1\}^{n_n}$ be a vector of binary decision variables, modelling the placement of a chlorine booster station at network nodes, i.e. $v^b_i=1$ if a chlorine booster station is placed at node $i$, $v^b_i=0$, otherwise. The number of boosters considered for installation is enforced by the linear constraint:
\begin{equation}
\label{eq:boost_con}
    \bm{1}^T\bm{v}^b = n_b 
\end{equation}
Chlorine concentration at unknown head node $i \in \mathcal{N}$ and time $k \in \mathcal{T}$ is governed by the following mixing equations:
\begin{subequations}
\label{eq:mixing_bigMs}
\begin{align}
    &c_{i,k}d_{i,k} + \sum_{l \in I^{\text{in}}_i}(w_{0,l,k}-\rho_{l,k}) + \sum_{l \in I^{\text{out}}_i}(\rho_{l,k}-w_{N_l,l,k}) - \xi_{i,k} = 0 \\
    &0 \leq \xi_{i,k} \leq \xi^{\max}_{i,k}v^b_i,
\end{align}
\end{subequations}
where slack variable $\xi_{i,k}\geq 0$ is {introduced} to model the additional constituent mass injected by a booster, and $\xi^{\max}_{i,k}$ are sufficiently large positive constants, for all $i \in \mathcal{N}$, and $k \in \mathcal{T}$. In addition, auxiliary variables $\rho_{l,k}$ are subject to the following linear constraints:
\begin{subequations}
 \label{eq:wq_auxvar}
 \begin{align}
& 0\leq \rho_{l,k} \leq \rho^{\max}_{l,k}z_{l,k} \\
&w_{0,l,k}+w_{N_l,l,k}-\rho_{l,k} - \rho^{\max}_{l,k}z_{l,k} \leq \rho^{\max}_{l,k}\\
&-w_{0,l,k}-w_{N_l,l,k}+\rho_{l,k}  \leq 0
\end{align}
\end{subequations}
Finally, we include the following lower and upper bounds:
\begin{subequations}
\label{eq:wq_bounds}
\begin{align}
&0 \leq c_{i,k} \leq c^{\max}_i, \quad i \in \mathcal{N} \cup \mathcal{N}^0, \; k \in \mathcal{T}, \\
& 0 \leq r_{j,l,k} \leq r^{\max}_{j,l}, \quad j=0,\ldots,J_l,\; l \in \mathcal{P},\; k \in \mathcal{T}, \\
& 0 \leq w_{j,l,k} \leq w^{\max}_{j,l,k} \quad j=0,\ldots,J_l,\; l \in \mathcal{P},\; k \in \mathcal{T},\\
& 0 \leq \rho_{l,k} \leq \rho^{\max}_{l,k}  \quad l \in \mathcal{P},\; k \in \mathcal{T},\\
&0 \leq \xi_{i,k} \leq \xi^{\max}_{i,k}, \quad i\in \mathcal{N},\; k \in \mathcal{T}.
\end{align}
\end{subequations}

{\subsection{Mixed Integer Non-linear Program}

The problem of optimal placement and control of valves and chlorine boosters aims to minimize \eqref{eq:objfun}, subject to non-convex quadratic constraints \eqref{eq:hyd_nonconvex} and \eqref{eq:wq_nonconvex}, and linear constraints \eqref{eq:ATD_auxvar}, \eqref{eq:binary_cons}, \eqref{eq:hyd_lincon}, \eqref{eq:hyd_aux_con}, \eqref{eq:hyd_bounds1},  \eqref{eq:hyd_bounds2}, \eqref{eq:upwind_disc}, \eqref{eq:boundary_disc}, \eqref{eq:boost_con}, \eqref{eq:mixing_bigMs}, \eqref{eq:wq_auxvar}, \eqref{eq:wq_bounds}. The optimization problem considers both continuous and binary variables. We write the problem in compact form, defining vectors $\bm{x}:=[\bm{q} \, \bm{h} \, \bm{\eta} \, \bm {\theta} ]^T$, $\bm{u}:=[\bm{s} \, \bm{q}^+ \, \bm{q}^- \, \bm{\theta}^+ \, \bm{\theta}^-]^T$, and $\bm{y}:=[\bm{c}\,\bm{r}\,\bm{w}\,\bm{\rho}\,\bm{\xi}\,\bm{\mu}]^T$. Consider the following Mixed Integer Non-linear Program (MINLP):
\begin{subequations}
\label{eq:prob_form}
\begin{alignat}{3}
&\underset{\substack{\bm{x},\bm{u},\bm{z}\\\bm{y},\bm{v},\bm{v}^b}}{\text{min}}&\; \; & f_{\text{AZP}}(\bm{x}) + f_{\text{ATD}}(\bm{y}) \label{eq:compact_obj}\\
&\text{s.t.}& & \bm{F}\bm{u} = \text{diag}(\bm{A}\bm{u})\bm{A}\bm{u} +\bm{B}\bm{u}\label{eq:compact_hyd_ncvx}\\
&&&\bm{W}\bm{y} = \text{diag}(\bm{S}\bm{u})\bm{R}\bm{y}\label{eq:compact_wq_bil}\\
&&& \bm{M}\bm{x} + \bm{N}\bm{u} + \bm{P}\bm{z} \leq \bm{p}\label{eq:compact_hyd_aux}\\
&&& \bm{x} \in X(\bm{v}) , \bm{v} \in V \label{eq:compact_valve_cons}\\
&&& \bm{y} \in Y(\bm{z},\bm{v}^b)\label{eq:compact_wq_cons}\\
&&& \bm{1}^T\bm{v}^b = n_b \\
&&& \bm{z} \in \{0,1\}^{n_tn_p}, \bm{v} \in \{0,1\}^{2n_p}, \bm{v}^b \in \{0,1\}^{n_n},
\end{alignat}
\end{subequations}
where, given a vector $\bm{e} \in \mathbb{R}^{N}$, $\text{diag}(\bm{e})\in \mathbb{R}^{N \times N}$ is the diagonal matrix with diagonal entries equal to the components of vector $\bm{e}$. Linear functions $f_{\text{AZP}}(\cdot)$ and $f_{\text{ATD}}(\cdot)$ are such that \eqref{eq:compact_obj} corresponds to \eqref{eq:objfun}. Matrices $\bm{F}$, $\bm{A}$, and $\bm{B}$ are defined so that the rows of \eqref{eq:compact_hyd_ncvx} correspond to the non-convex quadratic constraints \eqref{eq:hyd_nonconvex}. Matrices $\bm{W}$, $\bm{S}$, and $\bm{R}$ are opportunely defined so that the rows of \eqref{eq:compact_wq_bil} correspond to \eqref{eq:wq_nonconvex}. The set $V$ is defined by linear constraints \eqref{eq:binary_cons}. Given $\bm{v} \in V$, we denote by $X(\bm{v})$ the polyhedral set defined by constraints \eqref{eq:hyd_lincon} and \eqref{eq:hyd_bounds1}. Moreover, $\bm{M}$, $\bm{N}$, $\bm{P}$, and $\bm{p}$ are defined so that the rows of \eqref{eq:compact_hyd_aux} correspond to constraints \eqref{eq:hyd_aux_con} and \eqref{eq:hyd_bounds2}.
Finally, given vectors $\bm{v}^b$ and $\bm{z}$, we define $Y(\bm{v}^b,\bm{z})$ as the polyhedral set defined by constraints \eqref{eq:ATD_auxvar}, \eqref{eq:upwind_disc}, \eqref{eq:boundary_disc}, \eqref{eq:mixing_bigMs}, \eqref{eq:wq_auxvar}, \eqref{eq:wq_bounds}.
Problem \eqref{eq:prob_form} has $n_t(8n_p+4n_n+2\bar{J}+n_0)$ continuous variables, $n_tn_p + 2n_p+n_n$ binary variables, and $n_t(2n_p+\bar{J})$ non-convex quadratic constraints, where $\bar{J} = \sum_{l \in \mathcal{P}}(1+J_l)$. Therefore, even for small water networks, it results in large non-convex MINLPs, which are difficult to solve - see Table \ref{tab:prob_size}.}

{\section{Solution algorithm}
\label{sec:solalg}
The considered MINLP \eqref{eq:prob_form} combines difficulties in handling non-convex constraints with the presence of integer decision variables. In addition, the formulation of Problem \eqref{eq:prob_form} includes a discretized PDE for each network link, resulting in a large number of continuous variables and non-convex constraints, even for small size WDNs - see Table \ref{tab:prob_size}. 

We investigate the performance of off-the-shelf solvers to compute solutions for Problem \eqref{eq:prob_form}, considering two case study network models, namely \texttt{2loopsNet} and \texttt{pescara} - see Section \ref{sec:results} for network properties and layouts. We formulate Problem \eqref{eq:prob_form} in \texttt{2loopsNet} for $n_v \in \{1,2,3\}$ and $n_b=n_v$, and \texttt{pescara} for $n_v \in \{1,2,3,4,5\}$ and $n_b=n_v$ - a total of $8$ experiments. 
This study considers the global optimization solvers \texttt{BARON}~\cite{Tawarmalani2002}, \texttt{scip}~\cite{Gleixner2017The5.0}, \texttt{LINDOGlobal}~\cite{lindo}, and \texttt{Couenne}~\cite{Belotti2009}. Moreover, we investigate the ability of solvers \texttt{Bonmin}~\cite{Bonami2008}, \texttt{Knitro}~\cite{knitro}, \texttt{Ipopt}~\cite{Wachter2006}, and \texttt{AlphaECP}~\cite{alphaecp} to compute feasible solutions to Problem \eqref{eq:prob_form}. We refer to these as local solvers, because they do not provide guarantees of global optimality when considering non-convex MINLPs like Problem \eqref{eq:prob_form}. All experiments are performed on NEOS Server for Optimization \cite{czyzyk_et_al_1998}, with a time limit of $6$ hours. Since \texttt{Ipopt} does not directly handle problems with binary constraints, we have substituted them with the following complementary constraints:
\begin{equation}
    \begin{split}
        &\text{diag}(\bm{z})(\bm{1}-\bm{z}) = \bm{0}\\
        &\text{diag}(\bm{v})(\bm{1}-\bm{v}) = \bm{0}\\
        &\text{diag}(\bm{v^b})(\bm{1}-\bm{v^b}) = \bm{0}\\
    \end{split}
\end{equation}
The results of these experiments are summarized in Tables A1 - A16 of Appendix 2. When considering the small case study \texttt{2loopsNet}, \texttt{BARON} and \texttt{SCIP} were able to compute feasible solutions for $n_v \in \{2,3\}$. Moreover, local solvers \texttt{AlphaECP} and \texttt{Bonmin} have computed feasible solutions only when $n_v=2$. In comparison, \texttt{LINDOGlobal}, \texttt{Couenne}, \texttt{Knitro} and \texttt{Ipopt} failed to compute feasible solutions for all problem instances considering \texttt{2loopsNet}.
Furthermore, none of the tested solvers was able to compute feasible solutions for problem instances considering \texttt{pescara}. As off-the-shelf solvers were not able to compute feasible solutions in most tested instances, we propose an algorithm to compute feasible solutions for Problem \eqref{eq:prob_form}, together with bounds on their optimality gaps.

We propose the relax-tighten-round (RTR) algorithm, which combines a rounding heuristic with the solution of a continuous polyhedral relaxation of the non-convex MINLP in Problem~\eqref{eq:prob_form}, tightened using an optimization-based bound-tightening (OBBT) scheme. If successful, the algorithm computes a lower bound $\text{LB}$ and an upper bound $\text{UB}$ to the optimal value of Problem~\eqref{eq:prob_form}. A worst-case estimate on optimality gap of the computed solution is given by:
\begin{equation}
    \text{Gap}:= 100\frac{\text{UB} - \text{LB} }{\text{LB} }
\end{equation}
{In order to evaluate the obtained bounds on the optimality gaps, it is important to take into account the range of uncertainties that are inherent in hydraulic and water quality modelling of water networks. For example, \cite{Wright2015} and \cite{Waldron2020} showed that uncertainties affecting pressure control of operational water networks can result in up to $20 \%$ relative difference between simulated and measured pressure at network nodes. We expect the uncertainty range to be of the same magnitude and possibly higher for chlorine residuals.}

The steps necessary to derive the RTR algorithm are detailed in the following sub-sections. Section~\ref{sec:rounding_heuristic} presents a rounding heuristic to compute feasible solutions of Problem \eqref{eq:prob_form}. In Section~\ref{sec:pol_relax}, we introduce polyhedral relaxations of the non-convex constraints in Problem \eqref{eq:prob_form}. Then, Section~\ref{sec:obbt} describes the OBBT procedure to tighten the relaxation, and Section~\ref{sec:rtr_alg} presents the overall RTR algorithm.

\subsection{Rounding heuristic}
\label{sec:rounding_heuristic}
Firstly, we describe a heuristic to compute a feasible solution of Problem~\eqref{eq:prob_form}, given a vector of fractional values  $\bm{v} \in [0,1]^{2n_p} \cap V$. 
Let $J_{n_v} \subset \{1,\ldots,2n_p\}$ be the set of indices corresponding to the $n_v$ largest elements in $\bm{v}$, where only the largest value between ${v}_l$ and ${v}_{n_p+l}$ is considered for each link $l=1,\ldots,n_p$. For all $l=1,\ldots,2n_p$, define:
\begin{equation}
\label{eq:rounding}
\hat{v}_l = \begin{cases} 1 &\text{if } l \in J_{n_v} \\
0 &\text{otherwise}.
\end{cases}
\end{equation} 
This rounding scheme yields a vector $\bm{\hat{v}} \in V \cap \{0,1\}^{2n_p}$. Next, we obtain a feasible solution of Problem \eqref{eq:prob_form}. Observe that only constraints \eqref{eq:compact_wq_bil} and \eqref{eq:compact_wq_cons} couple vectors of hydraulic variables $\bm{x},\bm{u},\bm{z}$ with water quality vectors $\bm{y},\bm{v}^b$. We implement a two-stage approach, where hydraulic and water quality quantities are optimized in sequence. We consider the following MINLP:
\begin{subequations}
\label{eq:prob_AZP1}
\begin{alignat}{3}
&\underset{\substack{\bm{x},\bm{u},\bm{z}}}{\text{min}}&\; \; & f_{\text{AZP}}(\bm{x}) \label{eq:azp_obj1}\\
&\text{s.t.}& & \bm{F}\bm{u} = \text{diag}(\bm{A}\bm{u})\bm{A}\bm{u} +\bm{B}\bm{u}\label{eq:azp_hyd_ncvx1}\\
&&& \bm{M}\bm{x} + \bm{N}\bm{u} + \bm{P}\bm{z} \leq \bm{p}\label{eq:azp_hyd_aux1}\\
&&& \bm{x} \in X(\bm{\hat{v}}) , \bm{z} \in \{0,1\}^{n_tn_p}.
\end{alignat}
\end{subequations}
Problem \eqref{eq:prob_AZP1} includes $n_tn_p$ integer variables, and it is difficult to solve even for small/medium water networks. We have observed that $\bm{z} \in \{0,1\}^{n_tn_p}$, \eqref{eq:hyd_aux_con}, and \eqref{eq:hyd_nonconvex} are equivalent to non-linear equations \eqref{eq:hyd_nlp} and \eqref{eq:hyd_nlpa}. Since \eqref{eq:azp_hyd_ncvx1} and \eqref{eq:azp_hyd_aux1} correspond to constraints \eqref{eq:hyd_nonconvex} and \eqref{eq:hyd_aux_con}, respectively, Problem \eqref{eq:prob_AZP1} is equivalent to the following non-linear program:
\begin{subequations}
\label{eq:prob_AZP2}
\begin{alignat}{3}
&\underset{\substack{\bm{x}}}{\text{min}}&\; \; & f_{\text{AZP}}(\bm{x}) \label{eq:azp_obj2}\\
&\text{s.t.}& & \bm{g}(\bm{x})= \bm{0}\label{eq:azp_ncvx2}\\
&&& \bm{x} \in X(\bm{\hat{v}})  \label{eq:azp_valve_cons2}
\end{alignat}
\end{subequations}
where $\bm{g}(\cdot)$ is a non-linear function such that $\bm{g}(\bm{x})$ is the vector whose components are the rows of equalities \eqref{eq:hyd_nlpa}.
Let $\bm{\hat{x}} = [\bm{\hat{q}}\,\bm{\hat{h}}\,\bm{\hat{\eta}}\,\bm{\hat{\theta}}]^T$ be a locally optimal solution of Problem~\eqref{eq:prob_AZP2} computed by a NLP solver. We recover a feasible solution of Problem~\eqref{eq:prob_AZP1} by defining vectors $\bm{\hat{u}}=[\bm{\hat{s}}\,\bm{\hat{q}^+}\,\bm{\hat{q}^-}\,\bm{\hat{\theta}^+}\,\bm{\hat{\theta}^-}]^T$ and  $\bm{\hat{z}}$ using~\eqref{eq:hyd_nlp}. Finally, let $(\bm{\hat{y}},\bm{\hat{v}}^b)$ be solution of the mixed integer linear program (MILP):
\begin{equation}
\label{eq:prob_wqonly}
\begin{alignedat}{3}
&\underset{\substack{\bm{y},\bm{v}^b}}{\text{min}}&\; \; & f_{\text{ATD}}(\bm{y})\\
&\text{s.t.}& & \bm{W}\bm{y} = \text{diag}(\bm{S}\bm{\hat{u}})\bm{R}\bm{y}\\
&&& \bm{y} \in Y(\bm{\hat{z}},\bm{v}^b)\\
&&& \bm{1}^T\bm{v}^b = n_b \\
&&& \bm{v}^b \in \{0,1\}^{n_n}.
\end{alignedat}
\end{equation}
Since constraints in Problems \eqref{eq:prob_AZP1} and \eqref{eq:prob_wqonly} correspond to constraints in Problem~\eqref{eq:prob_form}, we conclude that $(\bm{\hat{x}},\bm{\hat{u}},\bm{\hat{z}},\bm{\hat{y}},\bm{\hat{v}},\bm{\hat{v}}^b$) is a feasible solution for Problem~\eqref{eq:prob_form}.
\subsection{Polyhedral relaxation}
\label{sec:pol_relax}
{In order to compute a lower bound to the optimal value of Problem \eqref{eq:prob_form}, we formulate a convex relaxation of the considered non-convex MINLP. Note that the non-convex terms in \eqref{eq:compact_hyd_ncvx} and \eqref{eq:compact_wq_bil} are the only non-linear terms within the formulation of Problem \eqref{eq:prob_form}. To take advantage of numerically efficient algorithms for solving linear programs, it is convenient to consider linear relaxations of these non-convex terms.} First, we consider polyhedral relaxations of matrix equation \eqref{eq:compact_hyd_ncvx}. Since each row in \eqref{eq:compact_hyd_ncvx} corresponds to a quadratic equation in \eqref{eq:hyd_nonconvex}, a polyhedral relaxation of \eqref{eq:compact_hyd_ncvx} is obtained by relaxing each row individually. As shown in Appendix 1, this results in the following linear equations:
\begin{subequations}
\label{eq:hyd_relax}
\begin{align}
    &\bm{F}\bm{u} \geq \text{diag}(\bm{A}\bm{u^{(i)}})\bm{A}(2\bm{u}-\bm{u}^{(i)})+\bm{B}\bm{u}, \; i=1,\ldots,m\label{eq:hyd_relax_a}\\
   &\bm{F}\bm{u} \leq \text{diag}(\bm{A}\bm{u}^{\min}+\bm{A}\bm{u}^{\max})\bm{A}\bm{u} + \bm{B}\bm{u} - \text{diag}(\bm{A}\bm{u}^{\min})\bm{A}\bm{u}^{\max}
    \end{align}
\end{subequations}
for given vectors $\bm{u}^{\min} = \bm{u}^{(1)} <\bm{u}^{(1)} <\ldots < \bm{u}^{(m)}= \bm{u}^{\max}$, where equality and inequality operators are to be interpreted element-wise. The polyhedral relaxations \eqref{eq:hyd_relax} have the advantage of resulting in linear programs, which can be efficiently solved by state-of-the art linear programming solvers.

The bilinear terms \eqref{eq:wq_nonconvex} are relaxed via the Reformulation Linearization Technique (RLT) \cite{Sherali1999a}. These relaxations are given by:
\begin{subequations}
\label{eq:wq_relax}
\begin{align}
   &\bm{W}\bm{y} \geq \text{diag}(\bm{S}\bm{u}^{\min})\bm{R}\bm{y} + \text{diag}(\bm{R}\bm{y}^{\min})\bm{S}\bm{u} - \text{diag}(\bm{S}\bm{u}^{\min})\bm{S}\bm{y}^{\min}\\
   &\bm{W}\bm{y} \geq \text{diag}(\bm{S}\bm{u}^{\max})\bm{R}\bm{y} + \text{diag}(\bm{R}\bm{y}^{\max})\bm{S}\bm{u} - \text{diag}(\bm{S}\bm{u}^{\max})\bm{S}\bm{y}^{\max}\\
   &\bm{W}\bm{y} \leq \text{diag}(\bm{S}\bm{u}^{\max})\bm{R}\bm{y} + \text{diag}(\bm{R}\bm{y}^{\min})\bm{S}\bm{u} - \text{diag}(\bm{S}\bm{u}^{\max})\bm{S}\bm{y}^{\min}\\
   &\bm{W}\bm{y} \leq \text{diag}(\bm{S}\bm{u}^{\min})\bm{R}\bm{y} + \text{diag}(\bm{R}\bm{y}^{\max})\bm{S}\bm{u} - \text{diag}(\bm{S}\bm{u}^{\min})\bm{S}\bm{y}^{\max}
\end{align}    
\end{subequations}
where inequalities are to be interpreted element-wise. Finally, observe that Problem~\eqref{eq:prob_form} includes a large number of binary variables, which results in impractical computational effort even for medium size water networks - see Table \ref{tab:prob_size}. Therefore, we consider the following continuous polyhedral relaxation of Problem~ \eqref{eq:prob_form}, where we also relax the binary constraints:
\begin{equation}
\label{eq:prob_relax}
\begin{alignedat}{3}
&\underset{\substack{\bm{x},\bm{u},\bm{z}\\\bm{y},\bm{v},\bm{v}^b}}{\text{min}}&\; \; & f_{\text{AZP}}(\bm{x}) + f_{\text{ATD}}(\bm{y}) \\
&\text{s.t.}& & \bm{F}\bm{u} \geq \text{diag}(\bm{A}\bm{u^{(i)}})\bm{A}(2\bm{u}-\bm{u}^{(i)})+\bm{B}\bm{u}, \; i=1,\ldots,m\\
&&&\bm{F}\bm{u} \leq \text{diag}(\bm{A}\bm{u}^{\min}+\bm{A}\bm{u}^{\max})\bm{A}\bm{u} + \bm{B}\bm{u} - \text{diag}(\bm{A}\bm{u}^{\min})\bm{A}\bm{u}^{\max}\\
&&&\bm{W}\bm{y} \geq \text{diag}(\bm{S}\bm{u}^{\min})\bm{R}\bm{y} + \text{diag}(\bm{R}\bm{y}^{\min})\bm{S}\bm{u} - \text{diag}(\bm{S}\bm{u}^{\min})\bm{S}\bm{y}^{\min}\\
&&&\bm{W}\bm{y} \geq \text{diag}(\bm{S}\bm{u}^{\max})\bm{R}\bm{y} + \text{diag}(\bm{R}\bm{y}^{\max})\bm{S}\bm{u} - \text{diag}(\bm{S}\bm{u}^{\max})\bm{S}\bm{y}^{\max}\\
&&&\bm{W}\bm{y} \leq \text{diag}(\bm{S}\bm{u}^{\max})\bm{R}\bm{y} + \text{diag}(\bm{R}\bm{y}^{\min})\bm{S}\bm{u} - \text{diag}(\bm{S}\bm{u}^{\max})\bm{S}\bm{y}^{\min}\\
&&&\bm{W}\bm{y} \leq \text{diag}(\bm{S}\bm{u}^{\min})\bm{R}\bm{y} + \text{diag}(\bm{R}\bm{y}^{\max})\bm{S}\bm{u} - \text{diag}(\bm{S}\bm{u}^{\min})\bm{S}\bm{y}^{\max}\\
&&& \bm{M}\bm{x} + \bm{N}\bm{u} + \bm{P}\bm{z} \leq \bm{p}\\\
&&& \bm{x} \in X(\bm{v}) , \bm{v} \in V \\
&&& \bm{y} \in Y(\bm{z},\bm{v}^b)\\
&&& \bm{1}^T\bm{v}^b = n_b \\
&&& \bm{z} \in [0,1]^{n_tn_p}, \bm{v} \in [0,1]^{2n_p}, \bm{v}^b \in [0,1]^{n_n}.
\end{alignedat}
\end{equation}
\subsection{{Optimization}-Based Bound-Tightening}
\label{sec:obbt}
The tightness of relaxations \eqref{eq:hyd_relax} and \eqref{eq:wq_relax}, depends on vectors $\bm{u}^{\min}$ and $\bm{u}^{\max}$, whose elements are lower and upper bounds on hydraulic auxiliary variables \eqref{eq:hyd_bounds2}. Equations \eqref{eq:hyd_nlp} and \eqref{eq:hyd_nlpa} imply that elements of $\bm{u}^{\min}$ and $\bm{u}^{\max}$ are functions of lower and upper bounds on the flow variables $q_{l,k}$, $l \in \mathcal{P}$, $k \in \mathcal{T}$. Hence, we consider an optimization-based bound-tightening (OBBT) scheme, to reduce the domain of the flow variables. We expect flow variables to be primarily influenced by hydraulic variables and constraints in Problem \eqref{eq:prob_form}. Hence. for all $\sigma \in \{-1,1\}$, $l \in \mathcal{P}$, and $k \in \mathcal{T}$, we consider: 
\begin{equation}
\label{eq:bt_large}
\begin{alignedat}{3}
&\underset{\substack{\bm{x},\bm{u},\bm{z},\bm{v}}}{\text{min}}&\; \; & \sigma \bm{e}^T_{(l,k)}\bm{x}\\
&\text{s.t.}& & \bm{F}\bm{u} \geq \text{diag}(\bm{A}\bm{u^{(i)}})\bm{A}(2\bm{u}-\bm{u}^{(i)})+\bm{B}\bm{u}, \; i=1,\ldots,m\\
&&&\bm{F}\bm{u} \leq \text{diag}(\bm{A}(\bm{u}^{\min}+\bm{u}^{\max}))\bm{A}\bm{u} + \bm{B}\bm{u}- \text{diag}(\bm{A}\bm{u}^{\min})\bm{A}\bm{u}^{\max}\\
&&& \bm{M}\bm{x} + \bm{N}\bm{u} + \bm{P}\bm{z} \leq \bm{p}\\\
&&& \bm{x} \in X(\bm{v}) , \bm{v} \in V \\
&&& \bm{z} \in [0,1]^{n_tn_p}, \bm{v} \in [0,1]^{2n_p}.
\end{alignedat}
\end{equation}
where $\bm{e}_{(l,k)}$ is an opportunely defined vector, which selects variable $q_{l,k}$ from vector $\bm{x}$. However, solving $2n_tn_p$ Problems \eqref{eq:bt_large} would require a significant computational effort even for small/medium size water networks. We investigate an alternative approach, aimed at solving smaller linear programs than \eqref{eq:bt_large}. Observe that matrices $\bm{F},\bm{A},\bm{B},\bm{M},\bm{N}$, and $\bm{P}$ are block diagonal with respect to the time index $k \in \mathcal{T}$, as they correspond to constraints \eqref{eq:hyd_aux_con} and \eqref{eq:hyd_nonconvex}. Let $\bm{F}_k,\bm{A}_k,\bm{B}_k,\bm{M}_k,\bm{N}_k$ and $\bm{P}_k$ be the diagonal blocks, respectively, for all $k \in \mathcal{T}$. In addition, let $\bm{x}_k,\bm{u}_k,\bm{z}_k$ be the sub-vectors of $\bm{x},\bm{u},\bm{z}$ whose components correspond to time index $k \in \mathcal{T}$.
We denote by $X_k(\bm{v})$ the polyhedral set defined by constraints \eqref{eq:hyd_lincon} considering only time index $k \in \mathcal{T}$. Finally, we introduce new vectors of variables $\bm{\tilde{v}}_k \in \mathbb{R}^{2n_p}$, $k \in \mathcal{T}$. The following problem is equivalent to Problem~\eqref{eq:bt_large}:
\begin{equation}
\label{eq:bt_1}
\begin{alignedat}{3}
&\underset{\substack{\bm{x},\bm{u},\bm{z}\\\bm{v},\bm{\tilde{v}}}}{\text{min}}&\; \; & \sigma \bm{e}^T_{(l,k)}\bm{x}\\
&\text{s.t.}& & \bm{F}_k\bm{u}_k \geq \text{diag}(\bm{A}_k\bm{u_k^{(i)}})\bm{A}_k(2\bm{u}_k-\bm{u}_k^{(i)})+\bm{B}_k\bm{u}_k, \; i=1,..,m, \; k \in \mathcal{T}\\
&&&\bm{F_k}\bm{u}_k \leq \text{diag}(\bm{A}_k(\bm{u}_k^{\min}+\bm{u}_k^{\max}))\bm{A}_k\bm{u}_k + \bm{B}_k\bm{u}_k - \text{diag}(\bm{A}_k\bm{u}_k^{\min})\bm{A}\bm{u}_k^{\max}, \; k \in \mathcal{T}\\
&&& \bm{M}_k\bm{x}_k + \bm{N}_k\bm{u}_k + \bm{P}_k\bm{z}_k \leq \bm{p}_k, \; k \in \mathcal{T}\\\
&&& \bm{x}_k \in X_k(\bm{\tilde{v}_k}), \bm{\tilde{v}}_k \in V, \; k \in \mathcal{T} \\
&&& \bm{z}_k \in [0,1]^{n_p},\bm{\tilde{v}}_k \in [0,1]^{2n_p}, \; k \in \mathcal{T}\\
&&& \bm{\tilde{v}}_k=\bm{v}, \; k \in \mathcal{T}.
\end{alignedat}
\end{equation}
Removing the time-coupling constraints on vectors $\bm{v}$ and $\bm{\tilde{v}}$, Problem~\eqref{eq:bt_1} becomes separable with respect to $k \in \mathcal{T}$. The considered OBBT scheme solves $2n_tn_p$ linear programs of the form:
\begin{equation}
\label{eq:bt_2}
\begin{alignedat}{3}
&\underset{\substack{\bm{x}_k,\bm{u}_k,\bm{z}_k,\bm{\tilde{v}}_k}}{\text{min}}&\; \; & \sigma \bm{e}^T_{(l)}\bm{x}_k\\
&\text{s.t.}& & \bm{F}_k\bm{u}_k \geq \text{diag}(\bm{A}_k\bm{u_k^{(i)}})\bm{A}_k(2\bm{u}_k-\bm{u}_k^{(i)})+\bm{B}_k\bm{u}_k, \; i=1,\ldots,m\\
&&&\bm{F_k}\bm{u}_k \leq \text{diag}(\bm{A}_k(\bm{u}_k^{\min}+\bm{u}_k^{\max}))\bm{A}_k\bm{u}_k + \bm{B}_k\bm{u}_k- \text{diag}(\bm{A}_k\bm{u}_k^{\min})\bm{A}\bm{u}_k^{\max}\\
&&& \bm{M}_k\bm{x}_k + \bm{N}_k\bm{u}_k + \bm{P}_k\bm{z}_k \leq \bm{p}_k\\
&&& \bm{x}_k \in X_k(\bm{\tilde{v}_k}), \bm{\tilde{v}}_k \in V\\ 
&&& \bm{z}_k \in [0,1]^{n_p}, \bm{\tilde{v}}_k \in [0,1]^{2n_p}.
\end{alignedat}
\end{equation}
where $\bm{e}_{(l)}$ is a vector used to select element $q_{l,k}$ from vector $\bm{x}_k$, for all $\sigma \in \{-1,1\}$, $k \in \mathcal{T}$, and $l \in \mathcal{P}$. The OBBT scheme is summarized in Algorithm \ref{alg:bt}.
\begin{algorithm}
\caption{Optimization-based bound-tightening (OBBT)}
\label{alg:bt}
\begin{algorithmic}[H]
\State \textbf{Output: } tightened vectors  $\bm{u}^{\min}$ and $\bm{u}^{\max}$.
\For{$l\in\mathcal{P}$}
\For{$k \in \mathcal{T}$}
\State Tighten $q_{l,k}^{\min}$ by solving Problem~\eqref{eq:bt_2} with $\sigma =1$.
\State Tighten $q_{l,k}^{\max}$ by solving Problem~\eqref{eq:bt_2} with $\sigma =-1$.
\State Update $\bm{u}^{\min}$ and $\bm{u}^{\max}$ using \eqref{eq:hyd_nlp} and \eqref{eq:hyd_nlpa}.
\EndFor
\EndFor
\end{algorithmic}
\end{algorithm}
\subsection{Algorithm implementation}
\label{sec:rtr_alg}
The RTR algorithm is summarized in Algorithm \ref{alg:rtr}. {At each iteration $i\in\{1,\ldots,I^{\max}\}$, we solve Problem~\eqref{eq:prob_relax}, computing the optimal value $\text{LB}(i)$, and a corresponding vector $\bm{v}^{(i)} \in V \cap [0,1]^{2n_p}$. Since Problem \eqref{eq:prob_relax} is a polyhedral relaxation of the original non-convex Problem \eqref{eq:prob_form}, we have that $\text{LB}(i)$ is a lower bound to the optimal value of Problem \eqref{eq:prob_form}.}
Next, we implement the rounding scheme to obtain a vector $\bm{\hat{v}}^{(i)} \in V\cap \{0,1\}^{2n_p}$, and apply a NLP solver to compute a locally optimal solution to Problem \eqref{eq:prob_AZP2} with $\bm{\hat{v}}=\bm{\hat{v}}^{(i)}$. If the NLP solver is successful, we store the computed locally optimal solution. If the termination criterion is not satisfied, we implement Algorithm \ref{alg:bt} and proceed with a new iteration. The iterative procedure stops either when the maximum number of iterations has been reached, or the relative change in lower bounds computed in consecutive iterations is smaller than the tolerance. Then, we select the locally optimal solution of Problem \eqref{eq:prob_AZP2} resulting in the smallest value of $f_{\text{AZP}}(\cdot)$, and refer to the corresponding vectors as $\bm{x}^*,\bm{u}^*,\bm{z}^*,\bm{v}^*$. We also select the largest lower bound value $\text{LB}$. Finally, we compute $\bm{y}^*,(\bm{v}^b)^*$ by solving the MILP \eqref{eq:prob_wqonly} with $\bm{\hat{u}} = \bm{u}^*$ and $\bm{\hat{z}}=\bm{z}^*$. When successful, the RTR algorithm terminates with a feasible solution for Problem \eqref{eq:prob_form}, given by $\bm{x}^*,\bm{u}^*,\bm{z}^*,\bm{v}^*,\bm{y}^*,(\bm{v}^b)^*$. {As consequence, we also obtain an upper bound to the optimal value of Problem \eqref{eq:prob_form}, given by $\text{UB}=f_{\text{AZP}}(\bm{x}^*)+f_{\text{ATD}}(\bm{y}^*)$. Observe that the NLP solver does not need to compute the globally optimal solution to the non-convex Problem \eqref{eq:prob_AZP2}, as it is sufficient to compute a feasible solution to generate an upper bound for Problem \eqref{eq:prob_form}.}

\begin{algorithm}
\caption{Relax-Tighten-Round (RTR)}
\label{alg:rtr}
\begin{algorithmic}[1]
\State{Initialize $\epsilon_{\text{tol}}$, $I_{\max}$, and set $i=1$.}
\For{$i=1,\ldots,I^{\max}$} 
\State{Solve Problem~\eqref{eq:prob_relax} computing $\text{LB}(i)$ (the optimal value) and $\bm{v}^{(i)}$.}
\State{Implement rounding scheme \eqref{eq:rounding} and obtain $\bm{\hat{v}}^{(i)}$.}
\State{Compute a locally optimal solution to Problem \eqref{eq:prob_AZP2} with $\bm{\hat{v}}=\bm{\hat{v}}^{(i)}$. }
\State {If successful, let $f^{(i)}_{\text{AZP}}$ be the corresponding objective function value.}
\If{$i\leq I^{\max}-1$ \textbf{ and } ($i=1$ \textbf{or} $\frac{|\text{LB}(i)-\text{LB}(i-1)|}{\text{LB}(i-1)} > \epsilon_{\text{tol}}$)}
\State Implement Algorithm \ref{alg:bt} to tighten $\bm{u}^{\min}$ and $\bm{u}^{\max}$.
\Else
\State{Terminate loop.}
\EndIf
\EndFor
\State{Set $f^*_{\text{AZP}} = \min_i f^{(i)}_{\text{AZP}}$ and $\text{LB} = \max_i \text{LB}(i)$.}
\State{Let $\bm{x}^*$ and $\bm{v}^*$ be the vectors corresponding to the best locally optimal solution found, and define $\bm{u}^*$ and $\bm{z}^*$ using \eqref{eq:hyd_nlp}.
}
\State{Solve Problem \eqref{eq:prob_wqonly} with $\bm{\hat{z}}=\bm{z}^*$,  $\bm{\hat{u}}=\bm{u}^*$, obtaining $\bm{y}^*$ and~$(\bm{v}^b)^*$}.
\State Set $\text{UB} = f_{\text{AZP}}(\bm{x}^*) +f_{\text{ATD}}(\bm{y}^*)$
\end{algorithmic}
\end{algorithm}
}
\section{Case studies and results}
\label{sec:results}
We evaluate the RTR algorithm on different benchmark water distribution network models, with varying size and level of connectivity. {All LPs and MILPs are solved using the state-of-the-art solver \texttt{GUROBI} (v9.0) \cite{GurobiOptimization2020}, while the nonlinear programs are solved using the solver for large-scale optimization \texttt{Ipopt} (v3.12.9) \cite{Wachter2006}. In the implementation of \texttt{Ipopt}, we supply gradients and Jacobians to the solver, in order to take advantage of the sparse structure of optimization problems in water networks.}

\begin{figure}[h!]
    \centering
    \subfloat[\label{fig:2loops_layout}\texttt{2loopsNet}]{\includegraphics[width=0.35\textwidth]{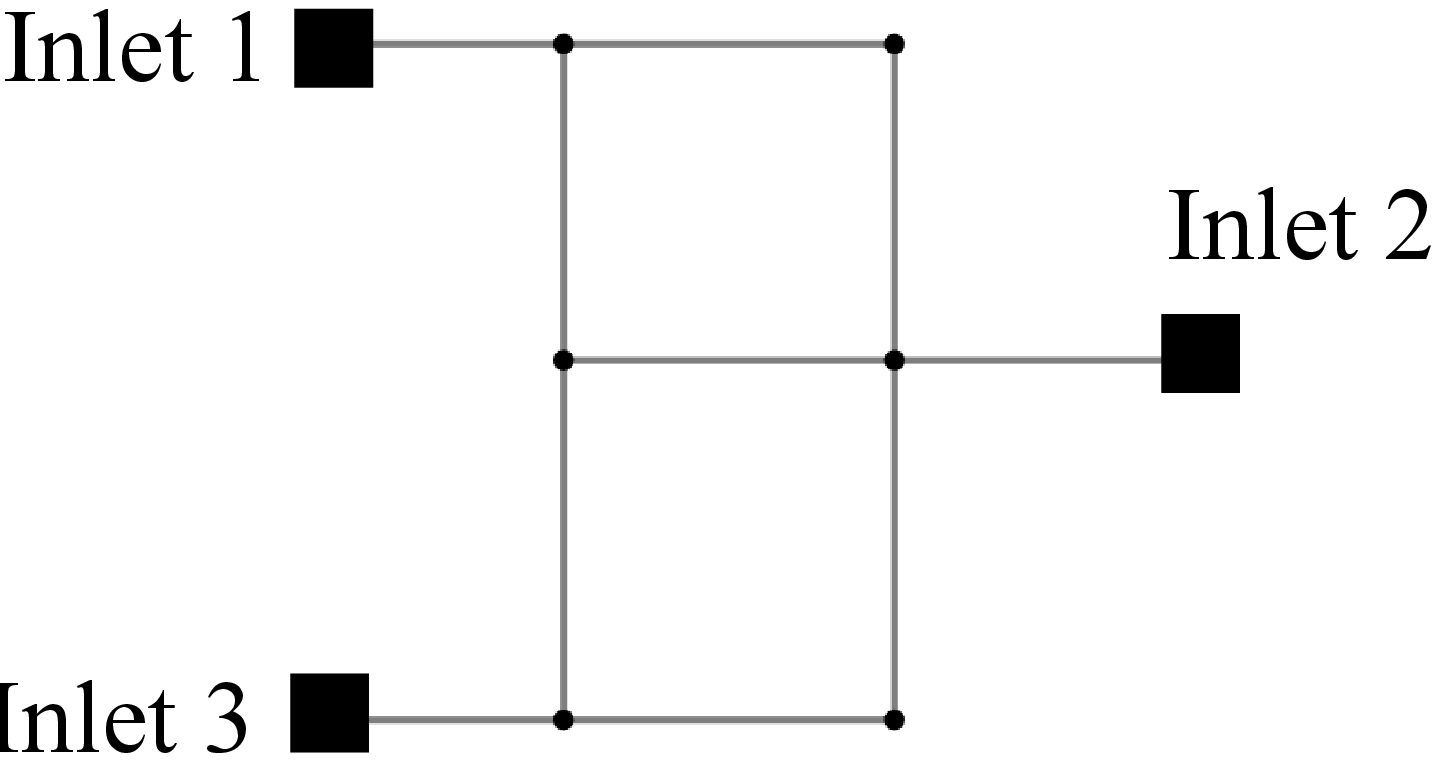}}\,
    \subfloat[\label{fig:pescara_layout}\texttt{pescara}]{\includegraphics[width=0.35\textwidth]{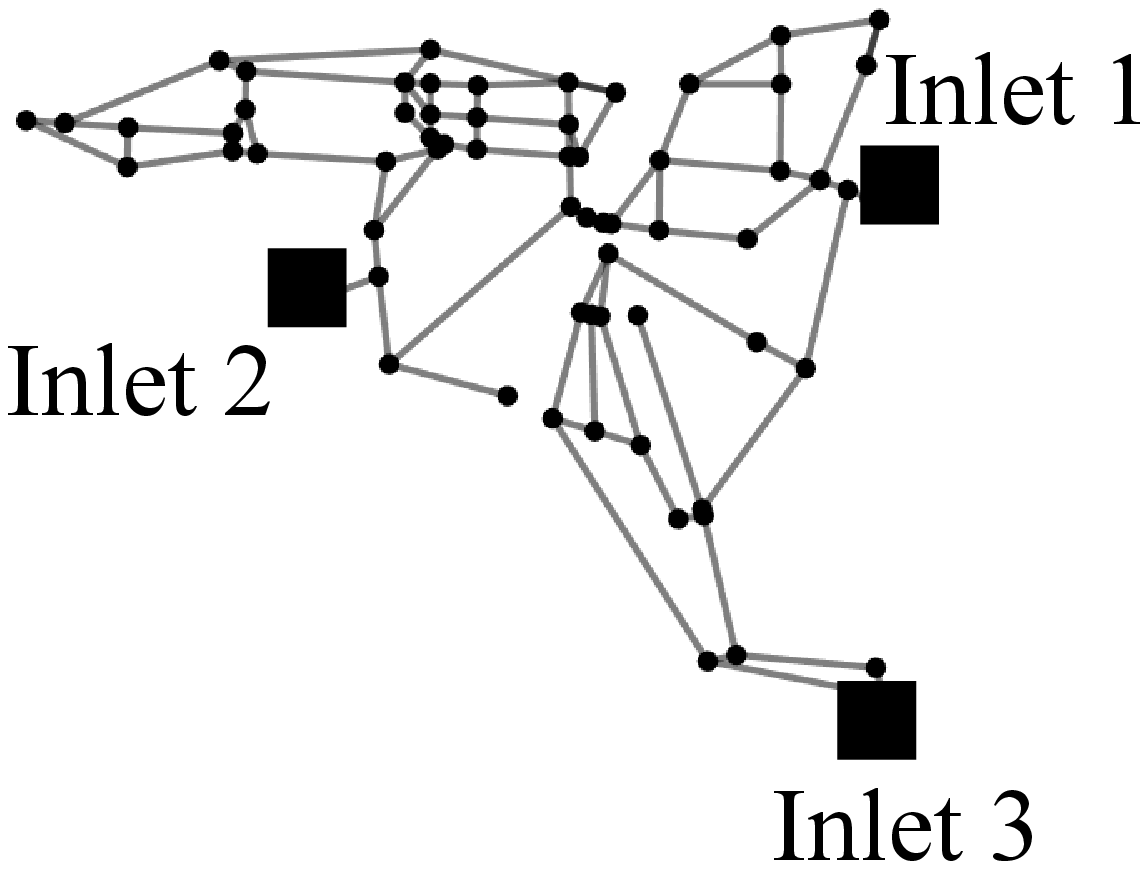}}\\
    \subfloat[\label{fig:modena_layout}\texttt{modena}]{\includegraphics[width=0.35\textwidth]{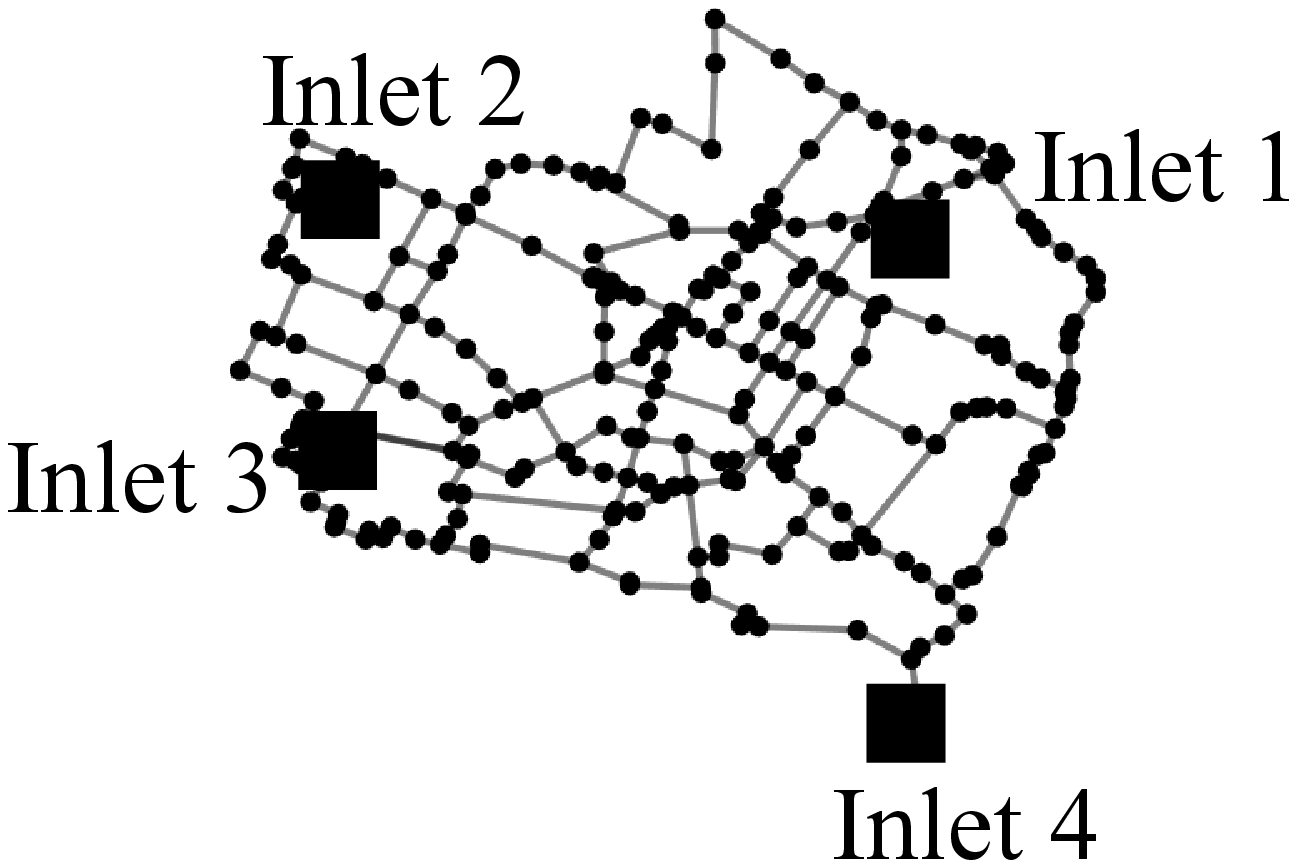}}\\    
    \caption{Case study network layouts.}
    \label{fig:networks}
\end{figure}

We consider a published benchmark network, referred to as \texttt{2loopsNet} \cite{Ostfeld2005}. In addition, we formulate and solve the problem of optimal valve and chlorine booster placement using \texttt{pescara} and \texttt{modena}, originally presented by \cite{Bragalli2011}. In order to obtain more realistic problem instances, we have introduced temporal and spatial variability of demand profiles and Hazen-Williams roughness coefficients, respectively. We have also added first-order chlorine decay coefficients to network pipes. All case study models consider $24$ hours of network operation, with a time step of one hour (i.e. $n_t=24$). Network hydraulic models and bounds on hydraulic heads and pipe flows are provided at \url{http://dx.doi.org/10.17632/ws9pwxkbb2.3}. The layout of \texttt{2loopsNet} is presented in Figure \ref{fig:2loops_layout}. The network has $n_n=6$ demand nodes, $n_p=10$ links, and $n_0=3$ water inlets. Case study \texttt{pescara} include $n_n=68$ nodes and $n_p=99$ links, and $n_0=3$ water inlets - see Figure \ref{fig:pescara_layout}. Finally, for \texttt{modena} we have $n_n=268$, $n_p=317$, and $n_0=4$ - see Figure \ref{fig:modena_layout}. {Note that} the considered case study networks result in large non-convex MINLPs, with a significant number of binary variables and non-convex terms - see Table \ref{tab:prob_size}.

In order to initialize nodal chlorine concentrations, we simulate $24$ hours of network operation using the software for hydraulic and water quality analysis EPANET~\cite{Rossman1996}, with fixed chlorine concentrations at inlets equal to $0.5$ mg/l. We define $c^0$ in equation \eqref{eq:upwind_disc} as the nodal concentrations at $24$:$00$ hours computed by the EPANET simulation. In the formulation of Problem~\eqref{eq:prob_form}, maximum allowed concentration at demand nodes is set $2$ mg/l, while chlorine concentrations at network inlets are not allowed to be greater than $0.5$ mg/l. The target concentration at demand nodes is $1$ mg/l. Finally, in \eqref{eq:upwind_disc}, we set a temporal time step $\Delta t=3600 s$ (1 hour) and $\Delta x_l=\frac{L_l}{2}$, where $L_l$ is the length of link $l$, for all $l=1,\ldots,n_p$. 
 { We set $m=5$ in \eqref{eq:hyd_relax} as preliminary experiments have shown that this setting results in sufficiently tight polyhedral relaxations, for all case studies. In addition, we have observed that the lower bounds computed by the RTR algorithm do not significantly improve after the first few iterations. Hence, we set $\epsilon_\text{tol}=10^{-2}$ and $I_{\max}=10$ in Algorithm \ref{alg:rtr}. Our choice is also supported by the results summarized in Table \ref{tab:computations}, which shows that the number of iterations performed by Algorithm \ref{alg:rtr} is never larger than $5$.}
\begin{table}
    \caption{Problem size for the $3$ case study networks.}
    \label{tab:prob_size}
    \centering
    \begin{tabular}{cccc}
    &$\#$ Cont. var. & $\#$ Bin. var. & $\#$ Non-convex cons. \\
    \hline
    \texttt{2loopsNet} &$4008$ & $266$ & $1200$ \\
    \texttt{pescara} & $39432$ &  $2615$ & $11760$ \\
    \texttt{modena} & $132336$ & $8510$ & $38040$
    \end{tabular}
\end{table}

 {
We formulate Problem \eqref{eq:prob_form} for $n_v \in \{1,2,3\}$ and $n_b \in \{0,\ldots,3\}$ in \texttt{2loopsNet}, and $n_v \in \{1,\ldots,5\}$ and $n_b \in \{0,\ldots,5\}$ in \texttt{pescara} and \texttt{modena}. Hence, we consider a total of $72$ different formulations of Problem \eqref{eq:prob_form}, and we implement the RTR Algorithm \ref{alg:rtr} to compute feasible solutions with bounds on their optimality gaps - see Appendix 2 in the supplementary material for tables of results.
In contrast to the off-the-shelf solvers considered in Section~\ref{sec:solalg}, RTR has computed feasible solutions in all problem instances for \texttt{2loopsNet} and \texttt{pescara}. In the case of \texttt{modena}, RTR has not returned a feasible solution only when $n_v=5$. 
{As reported in Appendix 2, the computed relative optimality gaps in problem instances for \texttt{2loopsNet} and \texttt{modena} are never larger than $20\%$. In comparison, in the case of \texttt{pescara}, the relative optimality gaps are between $25\%$ and $35\%$. Optimality gaps of such magnitude are comparable to the order of uncertainty affecting hydraulic models of operational water networks \cite{Wright2015,Waldron2020}. Hence, the RTR algorithm results in good quality solutions for the vast majority of problem instances.}

Table \ref{tab:computations} reports the computational effort required by RTR to compute feasible solutions for the tested problem instances. We set $N_{\text{AZP}}$ equal to the number of iterations required by \texttt{Ipopt} when computing a locally optimal solution for Problem \eqref{eq:prob_AZP2}. Moreover, we denote by $N_{\text{RTR}}$ the number of iterations of RTR Algorithm \ref{alg:rtr}. The number of calls of the OBBT Algorithm \ref{alg:bt} is then $N_{\text{RTR}}-1$, while the number of \texttt{Ipopt} calls within Algorithm \ref{alg:rtr} is equal to $N_{\text{RTR}}$. Table~\ref{tab:computations} shows that \texttt{Ipopt} is able to compute locally optimal solutions to Problem~\eqref{eq:prob_AZP2} within $30$ iterations for most problem instances. Moreover, these results show that the computational time required by RTR algorithm for \texttt{2loops} and \texttt{pescara} is significantly smaller than six hours (21600 s), the time limit set for the off-the-shelf solvers considered in Section~\ref{sec:solalg}.
\begin{table}[h!]

\small
    \caption{ (a) Minimum, (b) mean, and (c) maximum number of iterations and CPU time for the considered solvers and algorithms. $N_{\text{AZP}}=$number of iterations of IPOPT when applied to Problem~\eqref{eq:prob_AZP2}, $N_{\text{RTR}}=$number of iterations in RTR algorithm.}
    \label{tab:computations}
    \centering
    \begin{tabular}{ccc}
    \texttt{2loops} & \texttt{pescara} & \texttt{modena} \\
    \begin{tabular}{cccc}
    \hline
      & (a) & (b)  & (c)  \\
     \hline
      $N_{\text{RTR}}$ &   $3$ & $3$& $3$\\
      $N_{\text{AZP}}$ &  $12$ & $12$& $12$ \\
      Time (s) &   $5$ & $6$& $7$
    \end{tabular}
    &
        \begin{tabular}{cccc}
    \hline
      & (a) & (b) & (c)\\
     \hline
      &   $3$ & $4$& $5$\\
       &  $15$ & $21$& $23$ \\
       &  $210$ & $677$& $1387$
    \end{tabular}
    &
        \begin{tabular}{cccc}
    \hline
      & (a) & (b) & (c) \\
     \hline
       &   $2$ & $3$& $4$\\
       &  $15$ & $28$& $100$ \\
       &   $1045$ & $3027$ & $10132$
    \end{tabular}
    \end{tabular}
\end{table}
}
We also compare the lower bounds obtained by RTR with those computed by the off-the-shelf global optimization solvers \texttt{BARON}, \texttt{scip}, \texttt{Couenne}, and \texttt{LINDOGlobal} - see Tables A1-A8 in Appendix 2. Recall that we implemented these solvers for solving instances of Problem \eqref{eq:prob_form} formulated for $n_v\in\{1,2,3\}$ in \texttt{2loppsNet}, and $n_v \in \{1,2,3,4,5\}$ in \texttt{pescara}. Let $\text{LB}^{\text{OTS}}$ be the largest lower bound computed by the off-the-shelf solvers, for each experiment. We denote with $\text{LB}^{\text{RTR}}$ the lower bound computed by RTR for the same problem instances. As reported in Table \ref{tab:comp_lb}, the off-the-shelf global optimization solvers computed slightly better lower bounds in the case of \texttt{2loospNet}, with the largest difference roughly equal to $3.4 \%$. In comparison, in the case of \texttt{pescara}, the lower bounds computed by RTR are up to $18\%$ tighter than the best lower bounds obtained by the off-the-shelf solvers. We conclude that RTR has enabled the computation of lower bounds that are comparable to those obtained by off-the-shelf global optimization solvers after six hours of computations on the NEOS server.
\begin{table}[h!]
    \caption{Comparison between lower bounds computed by RTR and off-the-shelf-solvers.}
    \label{tab:comp_lb}
    \begin{tabular}{cccc|cccc}
        &&\texttt{2loopsNet}  & & &&\texttt{pescara}  & \\
        \hline
        $n_v$ & $n_b$ & $\text{LB}^{\text{RTR}}$ & $\text{LB}^{\text{OTS}}$  & $n_v$ & $n_b$ & $\text{LB}^{\text{RTR}}$ & $\text{LB}^{\text{OTS}}$ \\
        \hline
        1 & 1 & 99.63 & 102.05 & 1 & 1 & 34.60 & 28.30\\
        2 & 2 & 95.66 & 98.37 & 2 & 2 & 26.44 & 23.75\\
        3 & 3 & 93.13 & 96.43 & 3 & 3 & 22.28 & 21.82\\
          &   &   &   & 4 & 4 & 21.04 & 21.45\\
          &   &   &   & 5 & 5 & 20.32 & 21.01\\
        \end{tabular}
\end{table}

\begin{figure}[h!]
    \centering
    \includegraphics[width=0.5\textwidth]{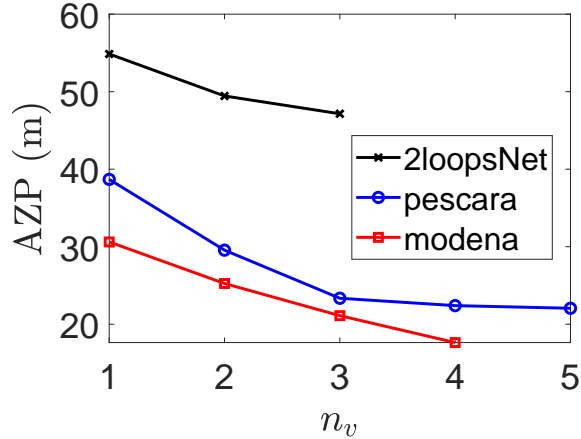}
\caption{AZP values for the considered problem instances.}
\label{fig:azp_opt}
\end{figure}

In Figure \ref{fig:azp_opt}, we report the computed AZP values for the feasible solutions obtained by RTR. For the same number of installed valves $n_v$, the computed AZP values for $n_b=0,\ldots,5$ are the same. As it should be expected, feasible solutions computed for increasing number of valves correspond to decreasing values of AZP.

\begin{figure}[h!]
    \centering
    \subfloat[\label{fig:ATD_2loops}\texttt{2loopsNet}]{
    \includegraphics[width=0.35\textwidth]{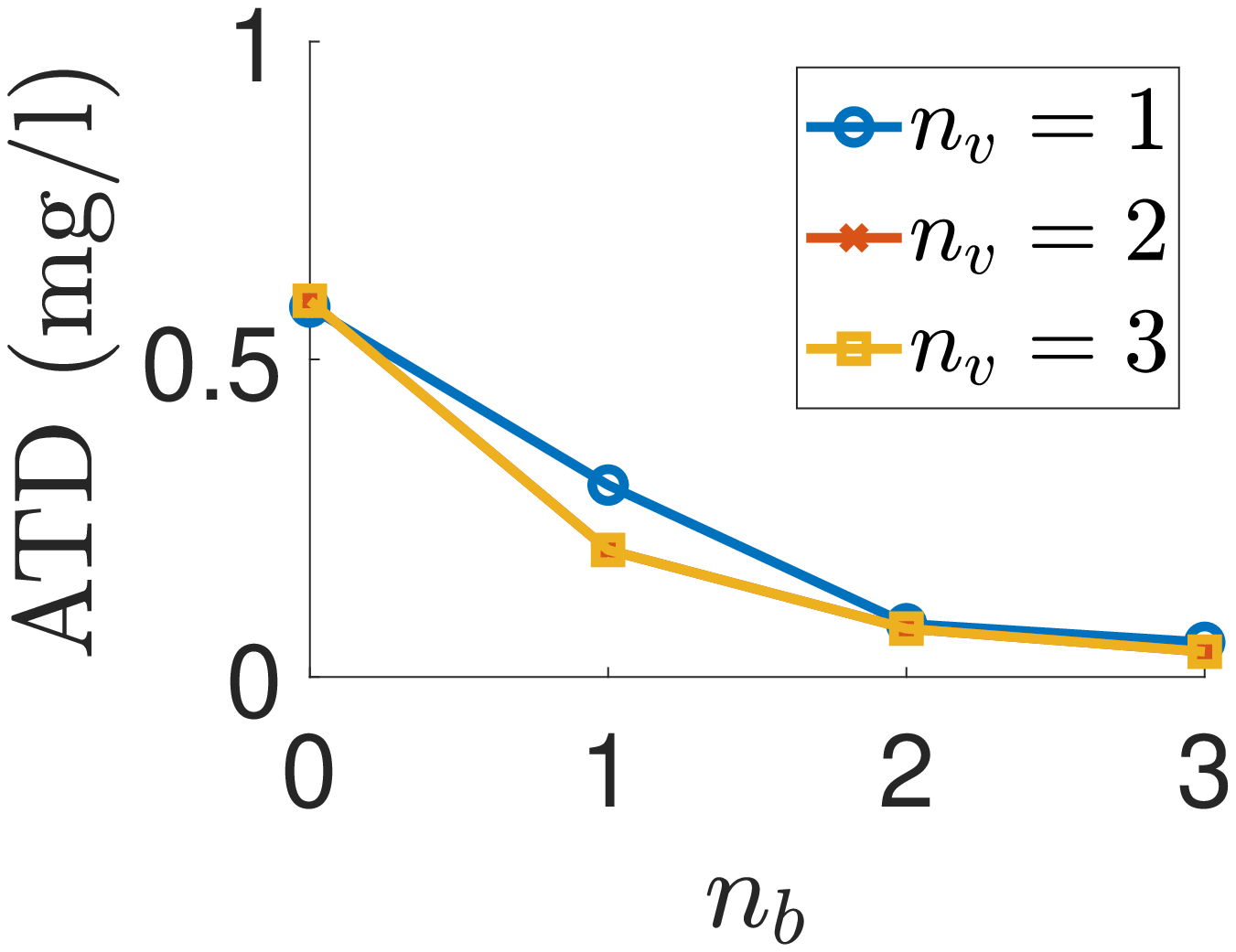}}
    \,
    \subfloat[\label{fig:ATD_pescara}\texttt{pescara}]{
    \includegraphics[width=0.35\textwidth]{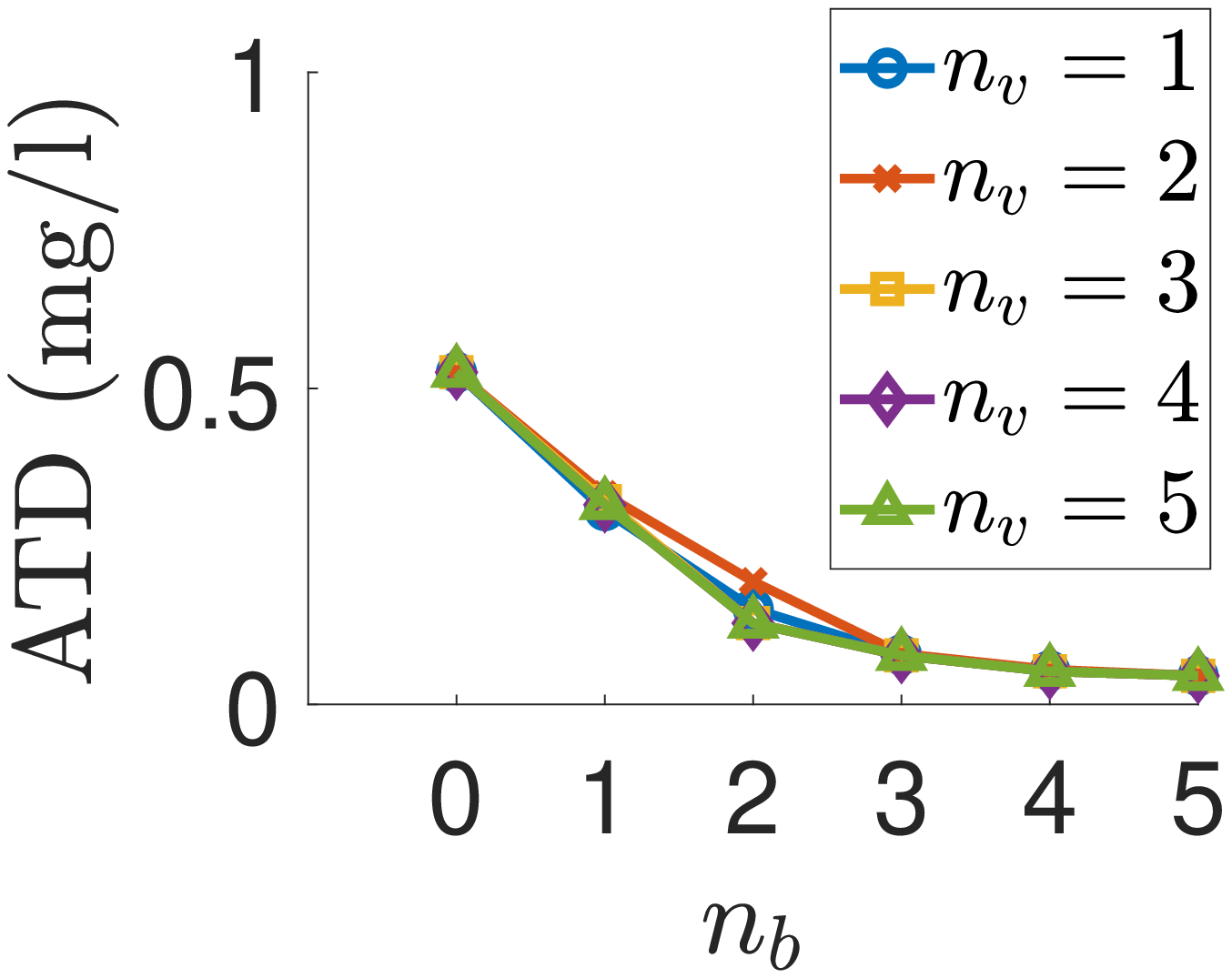}}
    \\
    \subfloat[\label{fig:ATD_modena}\texttt{modena}]{
    \includegraphics[width=0.35\textwidth]{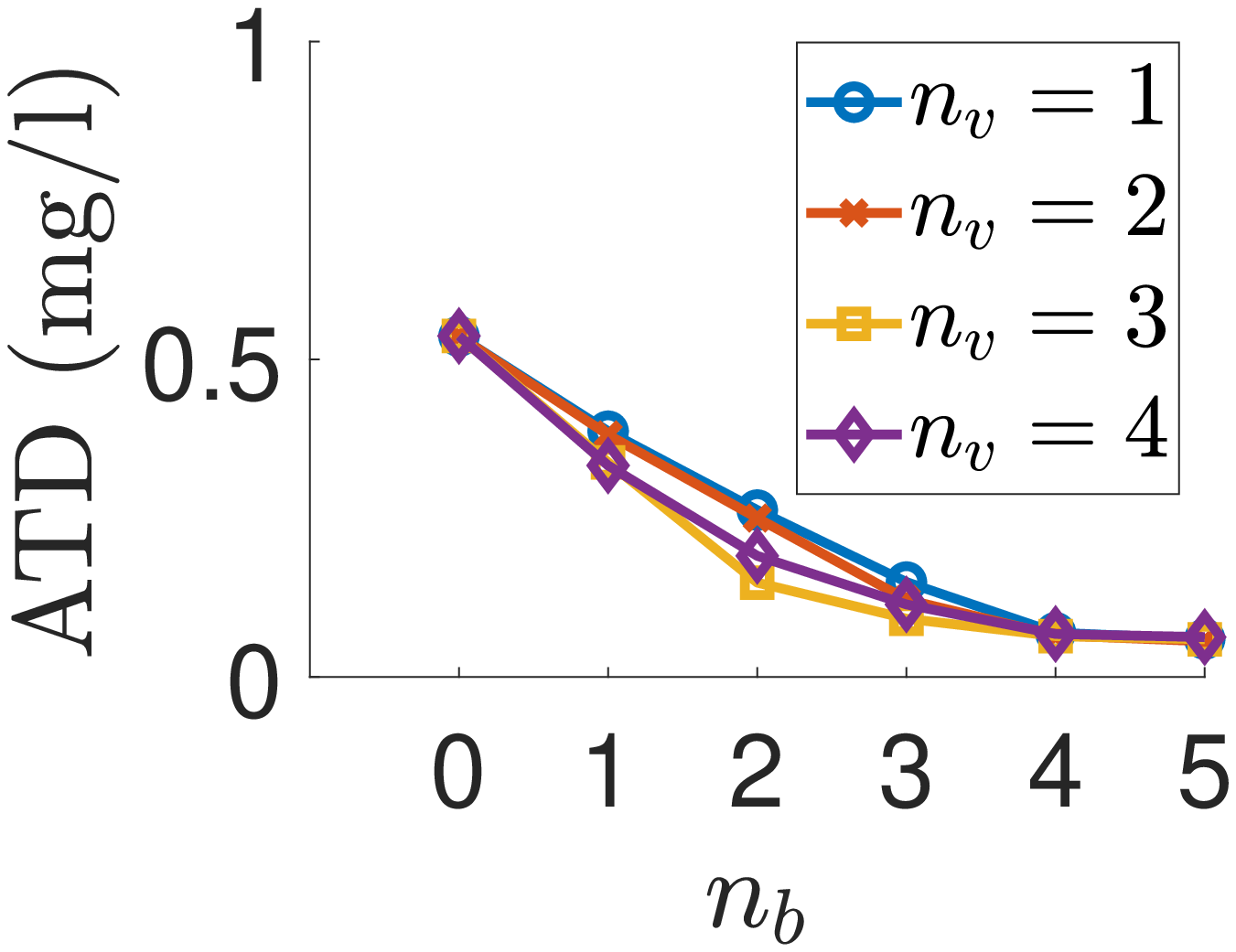}}
    \caption{optimized ATD values for the three case study networks.}
    \label{fig:atd_opt}
\end{figure}

Analogously, the average target deviation for nodal concentrations is reduced as additional chlorine booster stations are installed - see Figure \ref{fig:atd_opt}. Without any chlorine booster station, there is limited ability to control chlorine concentrations, using only the injected concentrations at water sources, which, in our formulation, can not be greater than $0.5$ (mg/l). As we install more booster stations, the system is able to maintain nodal concentrations closer to the target.

Observe that the best possible value of ATD is $0$. However, several hours are required for nodal concentrations to reach the optimized level, following the installation and operation of chlorine booster stations. If the travel time between a newly installed booster station and a specific node is $T^{\text{age}}$ hours, we expect nodal concentrations to reflect the action of the booster station after $T^{\text{age}}$ hours. In addition, network topology and spatial distribution of decay coefficients can affect the ability to control chlorine concentrations at selected locations. Therefore, we do not expect nodal concentrations to be exactly equal to the target at all time steps. 

\section{Conclusions}
We have proposed a new mixed integer nonlinear programming formulation for the problem of optimal placement and operation of pressure reducing valves and chlorine booster stations in water distribution networks. The numerical experiments reported in this manuscript show that off-the-shelf global optimization solvers can fail to compute feasible solutions for the considered problem, and the computed lower bounds to the optimal value are not tight. We have implemented polyhedral relaxations and a bound-tightening scheme resulting in improved lower bounds compared to off-the-shelf solvers. Furthermore, we have proposed the {Relax-Tighten-Round (RTR) algorithm} as heuristic to compute feasible solutions for the considered problem. The developed {RTR algorithm} has been evaluated by solving multiple problem instances for three case study networks. {RTR} is shown to outperform off-the-shelf solvers for the considered case studies. In addition, the developed heuristic has enabled the computation of good quality feasible solutions for the vast majority of the considered problem instances, {with bounds on the optimality gaps that are comparable to the order of uncertainty affecting hydraulic models of operational water networks.} 

{The proposed problem formulation and {RTR algorithm} enable the joint optimization of pressure and disinfectant dosage in water distribution networks. This allows water utilities to implement integrated and efficient schemes for pressure and water quality management, in order to minimize leakage and protect public health.
Future work should extend the problem formulation to include the operation of pumps and water tanks within the same optimization framework. Moreover, the proposed polyhedral relaxations could be tightened, for example implementing semidefinite or second-order cone relaxations of the non-convex quadratic constraints.}

\appendix
\renewcommand{\thetable}{A\arabic{table}}
\renewcommand{\thefigure}{A\arabic{figure}}
\renewcommand{\theequation}{A\arabic{equation}}
\section*{Appendix 1: polyhedral relaxation of quadratic head loss equation}
\begin{figure}
    \centering
    \subfloat[Convex quadratic relaxation \label{fig:cvx_relax}]{\includegraphics[width=0.4\textwidth]{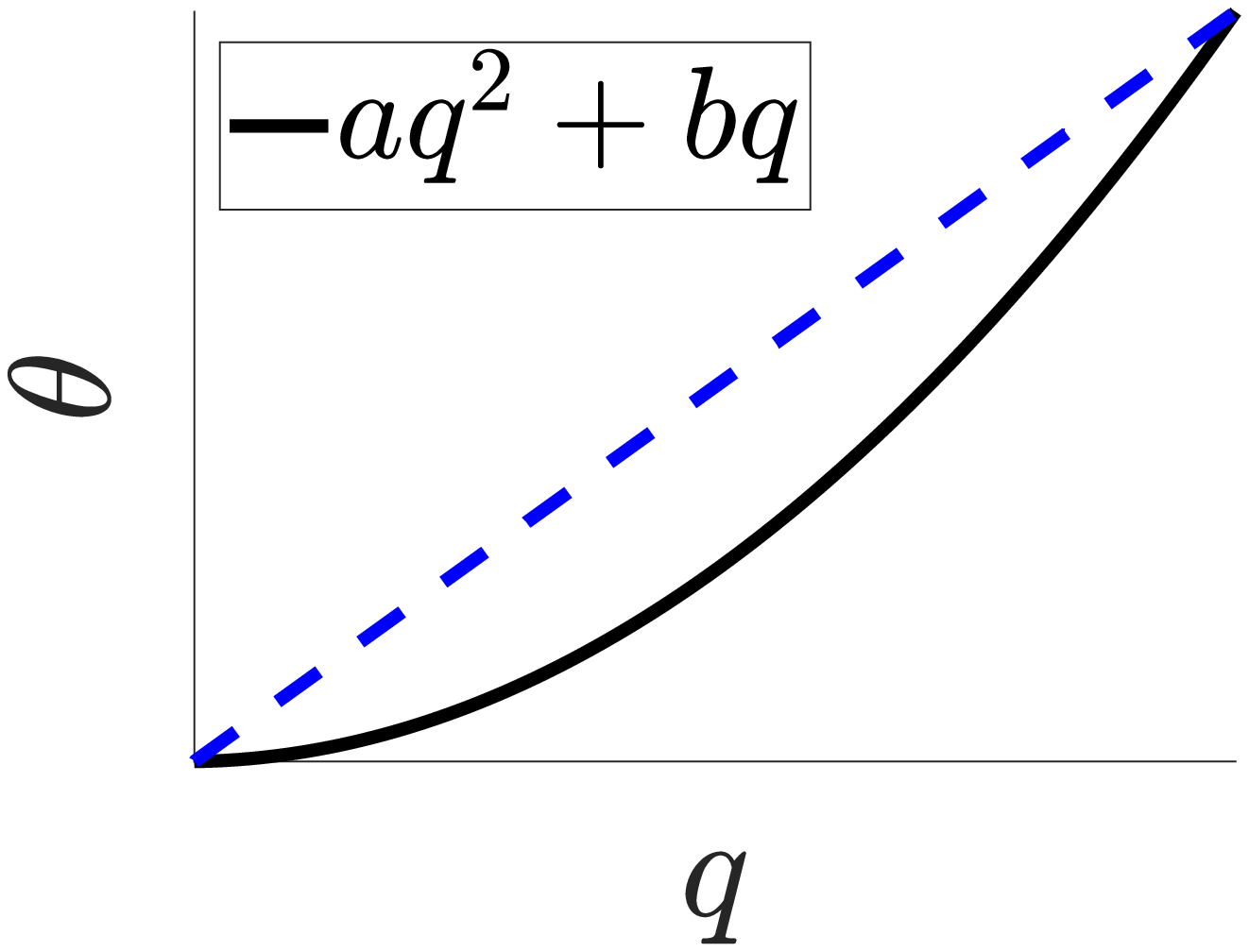}} \,
    \subfloat[Polyhedral relaxation \label{fig:lin_relax}]{\includegraphics[width=0.4\textwidth]{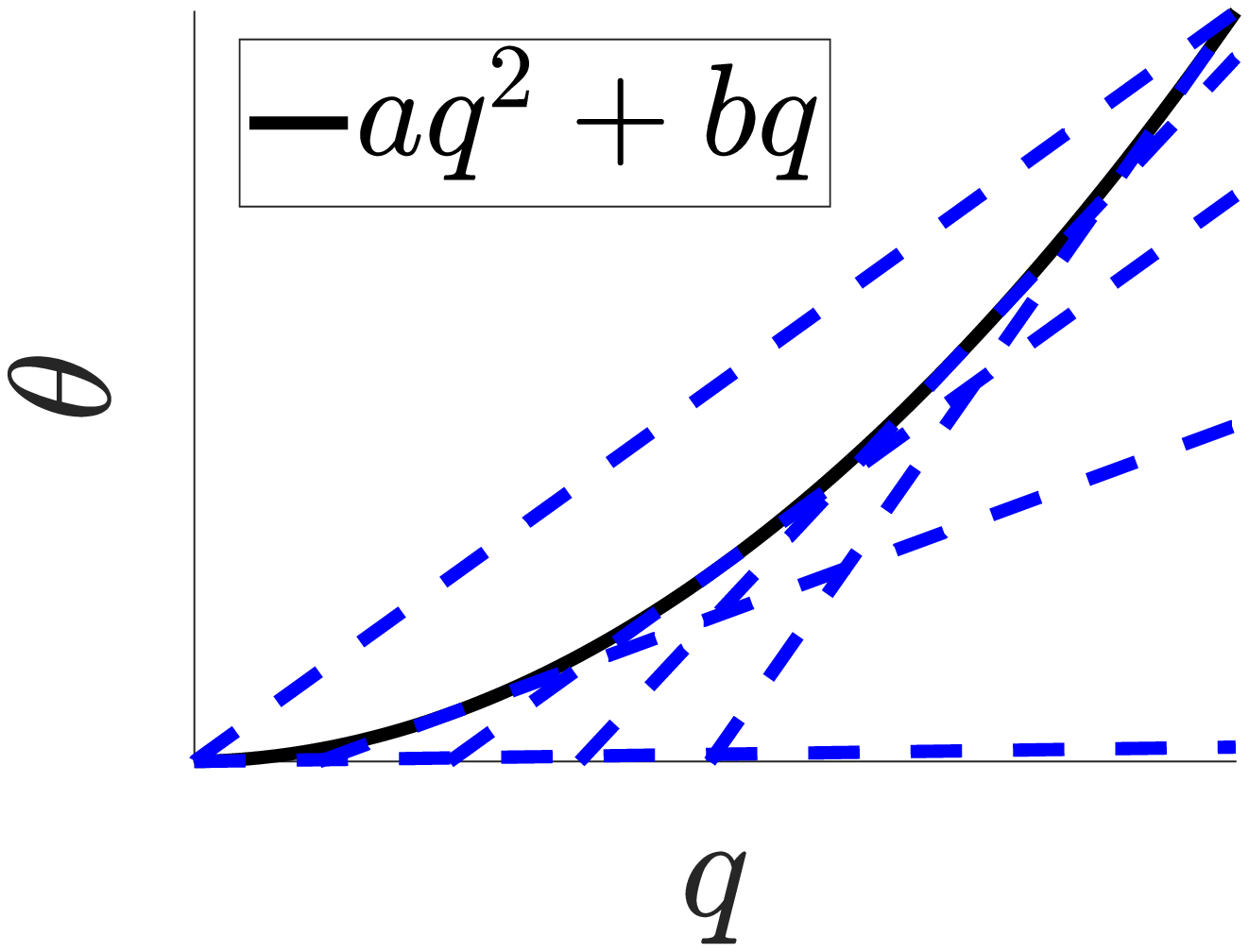}}
    \caption{Convex relaxations of the head loss equation.}
    \label{fig:my_label}
\end{figure}

We consider the non-convex quadratic constraint:
\begin{equation}
\label{eq:ncvx_quad}
    \theta = aq^2 + bq
\end{equation}
with $q \in [q^{\min},q^{\max}]$, $q^{\min}\geq 0$. A convex relaxation of \eqref{eq:ncvx_quad}, is given by:
\begin{equation}
\label{eq:cvx_relax}
    \begin{split}
    &\theta \geq aq^2 + bq \\
    &\theta \leq a(q^{\min}+q^{\max})q + bq- aq^{\min}q^{\max}
    \end{split}
\end{equation}
This is illustrated in Figure \ref{fig:cvx_relax}, where inequalities in \eqref{eq:cvx_relax} define the area between the curve and the dashed line. We can also further relax \eqref{eq:cvx_relax}, considering a linear outer approximation of the convex quadratic inequality constraint:
\begin{equation}
\label{eq:lin_relax}
    \begin{split}
    &\theta \geq a(q^{(i)})^2 + 2aq^{(i)}(q-q^{(i)}) + bq, \quad i=1,\ldots,m\\
    &\theta \leq a(q^{\min}+q^{\max})q + bq- aq^{\min}q^{\max}
    \end{split}
\end{equation}
where $q^{\min}=q^{(1)}<\ldots<q^{(m)}=q^{\max}$ are equidistant points with $m\geq2$. 
Equation \eqref{eq:lin_relax} defines a polyhedral relaxation of \eqref{eq:ncvx_quad} - Figure~\ref{fig:lin_relax} shows an example with $m=5$. 
Finally, observe that \eqref{eq:lin_relax} can equivalently re-written as:
\begin{equation}
\label{eq:lin_relax2}
    \begin{split}
    &\theta \geq a(q^{(i)})(2q-q^{(i)}) + bq, \quad i=1,\ldots,m\\
    &\theta \leq a(q^{\min}+q^{\max})q + bq - aq^{\min}q^{\max}
    \end{split}
\end{equation}
\section*{Appendix 2: tables of results}
\begin{table}[h!]
\caption{Results obtained by \texttt{BARON} for \texttt{2loopsNet}.}
\centering
\begin{tabular}{cccccc}
$n_v$ & $n_b$ & UB & LB & CPU Time (s) & Status \\ 
\hline 
1 & 1 & $-$ & 100.79 & 21600 & No solution \\ 
2 & 2 & 100.75 & 97.28 & 21600 & Integer solution \\ 
3 & 3 & 97.92 & 92.12 & 21600 & Integer solution \\ 
\hline 
\end{tabular}
\end{table}
\begin{table}[h!]
\caption{Results obtained by \texttt{BARON} for \texttt{pescara}.}
\centering
\begin{tabular}{cccccc}
$n_v$ & $n_b$ & UB & LB & CPU Time (s) & Status \\ 
\hline 
1 & 1 & $-$ & 22.56 & 21600 & No solution \\ 
2 & 2 & $-$ & 21.32 & 21600 & No solution \\ 
3 & 3 & $-$ & 20.65 & 21600 & No solution \\ 
4 & 4 & $-$ & 21.34 & 21600 & No solution \\ 
5 & 5 & $-$ & 21.01 & 21600 & No solution \\ 
\hline 
\end{tabular}
\end{table}
\begin{table}[h!]
\caption{Results obtained by \texttt{scip} for \texttt{2loopsNet}.}
\centering
\begin{tabular}{cccccc}
$n_v$ & $n_b$ & UB & LB & CPU Time (s) & Status \\ 
\hline 
1 & 1 & $-$ & 102.05 & 21600 & No solution \\ 
2 & 2 & 100.72 & 98.37 & 21600 & Integer solution \\ 
3 & 3 & 97.88 & 96.43 & 21600 & Integer solution \\ 
\hline 
\end{tabular}
\end{table}
\begin{table}[h!]
\caption{Results obtained by \texttt{scip} for \texttt{pescara}.}
\centering
\begin{tabular}{cccccc}
$n_v$ & $n_b$ & UB & LB & CPU Time (s) & Status \\ 
\hline 
1 & 1 & $-$ & 28.3 & 21600 & No solution \\ 
2 & 2 & $-$ & 23.76 & 21600 & No solution \\ 
3 & 3 & $-$ & 21.83 & 21600 & No solution \\ 
4 & 4 & $-$ & 21.46 & 21600 & No solution \\ 
5 & 5 & $-$ & 21.02 & 21600 & No solution \\ 
\hline 
\end{tabular}
\end{table}
\begin{table}[h!]
\caption{Results obtained by \texttt{LINDOGlobal} for \texttt{2loopsNet}.}
\centering
\begin{tabular}{cccccc}
$n_v$ & $n_b$ & UB & LB & CPU Time (s) & Status \\ 
\hline 
1 & 1 & $-$ & $-$ & 21600 & No solution \\ 
2 & 2 & $-$ & 95.4 & 21600 & No solution \\ 
3 & 3 & $-$ & 93.1 & 21600 & No solution \\ 
\hline 
\end{tabular}
\end{table}
\begin{table}[h!]
\caption{Results obtained by \texttt{LINDOGlobal} for \texttt{pescara}.}
\centering
\begin{tabular}{cccccc}
$n_v$ & $n_b$ & UB & LB & CPU Time (s) & Status \\ 
\hline 
1 & 1 & $-$ & $-$ & 21600 & No solution \\ 
2 & 2 & $-$ & $-$ & 21600 & No solution \\ 
3 & 3 & $-$ & $-$ & 21600 & No solution \\ 
4 & 4 & $-$ & $-$ & 21600 & No solution \\ 
5 & 5 & $-$ & $-$ & 21600 & No solution \\ 
\hline 
\end{tabular}
\end{table}
\begin{table}[h!]
\caption{Results obtained by \texttt{couenne} for \texttt{2loopsNet}.}
\centering
\begin{tabular}{cccccc}
$n_v$ & $n_b$ & UB & LB & CPU Time (s) & Status \\ 
\hline 
1 & 1 & $-$ & 95.01 & 21600 & No solution \\ 
2 & 2 & $-$ & 95 & 21600 & No solution \\ 
3 & 3 & $-$ & 92.71 & 21600 & No solution \\ 
\hline 
\end{tabular}
\end{table}
\begin{table}[h!]
\caption{Results obtained by \texttt{couenne} for \texttt{pescara}.}
\centering
\begin{tabular}{cccccc}
$n_v$ & $n_b$ & UB & LB & CPU Time (s) & Status \\ 
\hline 
1 & 1 & $-$ & 24.87 & 21600 & No solution \\ 
2 & 2 & $-$ & 22.12 & 21600 & No solution \\ 
3 & 3 & $-$ & 20.93 & 21600 & No solution \\ 
4 & 4 & $-$ & $-$ & 21600 & No solution \\ 
5 & 5 & $-$ & $-$ & 21600 & No solution \\ 
\hline 
\end{tabular}
\end{table}
\begin{table}[h!]
\caption{Results obtained by \texttt{bonmin} for \texttt{2loopsNet}.}
\centering
\begin{tabular}{ccccc}
$n_v$ & $n_b$ & UB & CPU Time (s) & Status \\ 
\hline 
1 & 1 & $-$ & 21600 & No solution \\ 
2 & 2 & 101.57 & 21600 & Integer solution \\ 
3 & 3 & $-$ & 21600 & No solution \\ 
\hline 
\end{tabular}
\end{table}
\begin{table}[h!]
\caption{Results obtained by \texttt{bonmin} for \texttt{pescara}.}
\centering
\begin{tabular}{ccccc}
$n_v$ & $n_b$ & UB & CPU Time (s) & Status \\ 
\hline 
1 & 1 & $-$ & 21600 & No solution \\ 
2 & 2 & $-$ & 21600 & No solution \\ 
3 & 3 & $-$ & 21600 & No solution \\ 
4 & 4 & $-$ & 21600 & No solution \\ 
5 & 5 & $-$ & 21600 & No solution \\ 
\hline 
\end{tabular}
\end{table}
\begin{table}[h!]
\caption{Results obtained by \texttt{Knitro} for \texttt{2loopsNet}.}
\centering
\begin{tabular}{ccccc}
$n_v$ & $n_b$ & UB & CPU Time (s) & Status \\ 
\hline 
1 & 1 & $-$ & 21600 & No solution \\ 
2 & 2 & $-$ & 21600 & No solution \\ 
3 & 3 & $-$ & 21600 & No solution \\ 
\hline 
\end{tabular}
\end{table}
\begin{table}[h!]
\caption{Results obtained by \texttt{Knitro} for \texttt{pescara}.}
\centering
\begin{tabular}{ccccc}
$n_v$ & $n_b$ & UB & CPU Time (s) & Status \\ 
\hline 
1 & 1 & $-$ & 21600 & No solution \\ 
2 & 2 & $-$ & 21600 & No solution \\ 
3 & 3 & $-$ & 21600 & No solution \\ 
4 & 4 & $-$ & 21600 & No solution \\ 
5 & 5 & $-$ & 21600 & No solution \\ 
\hline 
\end{tabular}
\end{table}
\begin{table}[h!]
\caption{Results obtained by \texttt{Ipopt} for \texttt{2loopsNet}.}
\centering
\begin{tabular}{ccccc}
$n_v$ & $n_b$ & UB & CPU Time (s) & Status \\ 
\hline 
1 & 1 & $-$ & 9.07 & Conv. local infeas. \\ 
2 & 2 & $-$ & 24.4 & Conv. local infeas. \\ 
3 & 3 & $-$ & 10.08 & Conv. local infeas. \\ 
\hline 
\end{tabular}
\end{table}
\begin{table}[h!]
\caption{Results obtained by \texttt{Ipopt} for \texttt{pescara}.}
\centering
\begin{tabular}{ccccc}
$n_v$ & $n_b$ & UB & CPU Time (s) & Status \\ 
\hline 
1 & 1 & $-$ & 2308.81 & Conv. local infeas. \\ 
2 & 2 & $-$ & 819.85 & Conv. local infeas. \\ 
3 & 3 & $-$ & 1592.12 & Conv. local infeas. \\ 
4 & 4 & $-$ & 1577.6 & Conv. local infeas. \\ 
5 & 5 & $-$ & 812.55 & Conv. local infeas. \\ 
\hline 
\end{tabular}
\end{table}
\begin{table}[h!]
\caption{Results obtained by \texttt{AlphaECP} for \texttt{2loopsNet}.}
\centering
\begin{tabular}{ccccc}
$n_v$ & $n_b$ & UB & CPU Time (s) & Status \\ 
\hline 
1 & 1 & $-$ & 21600 & No solution \\ 
2 & 2 & 100.19 & 21600 & Integer solution \\ 
3 & 3 & 97.9 & 21600 & Integer solution \\ 
\hline 
\end{tabular}
\end{table}
\begin{table}[h!]
\caption{Results obtained by \texttt{AlphaECP} for \texttt{pescara}.}
\centering
\begin{tabular}{ccccc}
$n_v$ & $n_b$ & UB & CPU Time (s) & Status \\ 
\hline 
1 & 1 & $-$ & 21600 & No solution \\ 
2 & 2 & $-$ & 21600 & No solution \\ 
3 & 3 & $-$ & 21600 & No solution \\ 
4 & 4 & $-$ & 21600 & No solution \\ 
5 & 5 & $-$ & 21600 & No solution \\ 
\hline 
\end{tabular}
\end{table}
\begin{table}[h!]
\caption{Results obtained by RTR for \texttt{2loopsNet}.}
\label{tab:results_2loops}
\centering
\begin{tabular}{cccccc}
$n_v$ & $n_b$ & Gap ($\%$) & $\text{UB}$ & $\text{LB}$ & CPU Time (s) \\ 
\hline 
1 & 0 & 6.26 & 106.05 & 99.8 & 6.68 \\ 
1 & 1 & 6.17 & 105.77 & 99.63 & 6.54 \\ 
1 & 2 & 5.95 & 105.56 & 99.63 & 6.29 \\ 
1 & 3 & 5.92 & 105.53 & 99.63 & 6.31 \\ 
2 & 0 & 5.08 & 100.65 & 95.79 & 5.89 \\ 
2 & 1 & 4.8 & 100.26 & 95.66 & 5.89 \\ 
2 & 2 & 4.67 & 100.13 & 95.66 & 5.94 \\ 
2 & 3 & 4.64 & 100.1 & 95.66 & 5.7 \\ 
3 & 0 & 5.48 & 98.35 & 93.24 & 5.48 \\ 
3 & 1 & 5.18 & 97.96 & 93.14 & 5.45 \\ 
3 & 2 & 5.05 & 97.84 & 93.14 & 5.5 \\ 
3 & 3 & 5.01 & 97.8 & 93.14 & 5.53 \\ 
\hline 
\end{tabular}
\end{table}
\begin{table}[h!]
\caption{Results obtained by RTR for \texttt{pescara}.}
\label{tab:results_pescara}
\centering
\small
\begin{tabular}{cccccc}
$n_v$ & $n_b$ & Gap ($\%$) & $\text{UB}$ & $\text{LB}$ & CPU Time (s) \\ 
\hline 
1 & 0 & 27.39 & 44.31 & 34.78 & 1304.52 \\ 
1 & 1 & 27.42 & 44.09 & 34.6 & 1324.84 \\ 
1 & 2 & 26.97 & 43.94 & 34.6 & 1310.12 \\ 
1 & 3 & 26.76 & 43.86 & 34.6 & 1311.94 \\ 
1 & 4 & 26.69 & 43.84 & 34.6 & 1314.38 \\ 
1 & 5 & 26.66 & 43.83 & 34.6 & 1386.77 \\ 
2 & 0 & 32.51 & 35.18 & 26.55 & 929.12 \\ 
2 & 1 & 32.32 & 34.99 & 26.44 & 953.04 \\ 
2 & 2 & 31.8 & 34.85 & 26.44 & 941.77 \\ 
2 & 3 & 31.37 & 34.74 & 26.44 & 945.98 \\ 
2 & 4 & 31.27 & 34.71 & 26.44 & 952.27 \\ 
2 & 5 & 31.24 & 34.7 & 26.44 & 1057.39 \\ 
3 & 0 & 29.32 & 28.98 & 22.41 & 588.54 \\ 
3 & 1 & 29.1 & 28.77 & 22.29 & 607.26 \\ 
3 & 2 & 28.23 & 28.58 & 22.29 & 598.21 \\ 
3 & 3 & 28 & 28.53 & 22.29 & 602.77 \\ 
3 & 4 & 27.89 & 28.5 & 22.29 & 604.91 \\ 
3 & 5 & 27.86 & 28.5 & 22.29 & 687.47 \\ 
4 & 0 & 32.59 & 28.03 & 21.14 & 216 \\ 
4 & 1 & 32.17 & 27.82 & 21.05 & 236.77 \\ 
4 & 2 & 31.28 & 27.63 & 21.05 & 230.53 \\ 
4 & 3 & 31.03 & 27.58 & 21.05 & 232.09 \\ 
4 & 4 & 30.91 & 27.55 & 21.05 & 230.9 \\ 
4 & 5 & 30.88 & 27.55 & 21.05 & 328.02 \\ 
5 & 0 & 35.65 & 27.68 & 20.41 & 210.3 \\ 
5 & 1 & 35.15 & 27.47 & 20.33 & 233.78 \\ 
5 & 2 & 34.23 & 27.29 & 20.33 & 223.26 \\ 
5 & 3 & 33.98 & 27.24 & 20.33 & 223.86 \\ 
5 & 4 & 33.85 & 27.21 & 20.33 & 224.49 \\ 
5 & 5 & 33.82 & 27.2 & 20.33 & 309.78 \\ 
\hline 
\end{tabular}
\end{table}
\begin{table}[h!]
\caption{Results obtained by RTR for \texttt{modena}.}
\label{tab:results_modena}
\small
\centering
\begin{tabular}{cccccc}
$n_v$ & $n_b$ & Gap ($\%$) & $\text{UB}$ & $\text{LB}$ & CPU Time (s) \\ 
\hline 
1 & 0 & 19.55 & 66.58 & 55.7 & 3961.52 \\ 
1 & 1 & 19.35 & 66.43 & 55.66 & 4054.75 \\ 
1 & 2 & 19.13 & 66.31 & 55.66 & 4197.41 \\ 
1 & 3 & 18.92 & 66.19 & 55.66 & 5005.25 \\ 
1 & 4 & 18.78 & 66.11 & 55.66 & 5064.04 \\ 
1 & 5 & 18.76 & 66.1 & 55.66 & 9002.53 \\ 
2 & 0 & 17.45 & 61.22 & 52.12 & 2382.6 \\ 
2 & 1 & 17.16 & 61.07 & 52.12 & 2488.16 \\ 
2 & 2 & 16.91 & 60.93 & 52.12 & 2812.11 \\ 
2 & 3 & 16.66 & 60.81 & 52.12 & 3275.83 \\ 
2 & 4 & 16.55 & 60.75 & 52.12 & 4287.42 \\ 
2 & 5 & 16.53 & 60.74 & 52.12 & 10131.5 \\ 
3 & 0 & 12.2 & 57.07 & 50.86 & 1106.79 \\ 
3 & 1 & 11.81 & 56.87 & 50.86 & 1229.86 \\ 
3 & 2 & 11.44 & 56.68 & 50.86 & 1223.84 \\ 
3 & 3 & 11.33 & 56.62 & 50.86 & 1469.3 \\ 
3 & 4 & 11.27 & 56.6 & 50.86 & 2273.58 \\ 
3 & 5 & 11.26 & 56.59 & 50.86 & 9419.76 \\ 
4 & 0 & 6.09 & 53.6 & 50.52 & 1065.34 \\ 
4 & 1 & 5.69 & 53.39 & 50.52 & 1135.22 \\ 
4 & 2 & 5.41 & 53.25 & 50.52 & 1189.47 \\ 
4 & 3 & 5.26 & 53.18 & 50.52 & 1352.26 \\ 
4 & 4 & 5.16 & 53.13 & 50.52 & 1794.46 \\ 
4 & 5 & 5.15 & 53.12 & 50.52 & 4573.93 \\ 
5 & 0 & $-$ & $-$ & 50.45 & 1048.34 \\ 
5 & 1 & $-$ & $-$ & 50.45 & 1045.1 \\ 
5 & 2 & $-$ & $-$ & 50.45 & 1052.41 \\ 
5 & 3 & $-$ & $-$ & 50.45 & 1051.53 \\ 
5 & 4 & $-$ & $-$ & 50.45 & 1050.66 \\ 
5 & 5 & $-$ & $-$ & 50.45 & 1049.48 \\ 
\hline 
\end{tabular}
\end{table}
\cleardoublepage
\section*{Acknowledgements}
\noindent Filippo Pecci and Ivan Stoianov are supported by EPSRC (EP/P004229/1, Dynamically Adaptive and Resilient Water Supply Networks for a Sustainable Future). Avi Ostfeld is supported by the Israel Science Foundation (grant No. 555/18).
\bibliography{references}

\end{document}